# THE EQUIVALENCE AND THE EMBEDDING PROBLEMS FOR CR-STRUCTURES OF ANY CODIMENSION

SUNG-YEON KIM AND DMITRI ZAITSEV

ABSTRACT. We give a solution to the equivalence and the embedding problems for smooth CR-submanifolds of complex spaces (and, more generally, for abstract CR-manifolds) in terms of complete differential systems in jet bundles satisfied by all CR-equivalences or CR-embeddings respectively (local and global). For the equivalence problem, the manifolds are assumed to be of finite type and finitely nondegenerate. These are higher order generalizations of the corresponding nondegeneracy conditions for the Levi form. It is shown by a simple example that these nondegeneracy conditions cannot be even slightly relaxed to more general known conditions. In particular, for essentially finite hypersurfaces in $\mathbb{C}^2$, such a complete system may not exist in general. For the embedding problem, the source manifold is assumed to be of finite type and the embeddings to be finitely nondegenerate. Situations are given, where the last condition is automatically satisfied by all CR-embeddings.

## CONTENTS

The first author was supported by postdoctoral fellowship program from Korea Science and Engineering Foundation (KOSEF).

## 1. INTRODUCTION

1.1. **The scope of the paper.** One of the main result of this paper is to give a solution to the equivalence problem for ($\mathcal{C}^\kappa$-smooth) CR-structures having nondegenerate Levi forms or, more generally, satisfying nondegeneracy conditions in terms of "higher order Levi forms". The *equivalence problem*, that goes back to H. POINCARÉ [P07] and E. CARTAN [CaE32], can be stated as follows:

*Given two manifolds with CR-structures, when does there exist a CR-diffeomorphism between them?*

This problem is an important case of the more general equivalence problem for the so-called *geometric structures* (see e.g. [St83]) introduced by E. CARTAN. It is also closely related to the biholomorphic equivalence problem for bounded domains in $\mathbb{C}^N$. Indeed, by FEFFERMAN's theorem [Fe74], two strongly pseudoconvex smoothly bounded domains are biholomorphically equivalent if and only if their boundaries have equivalent CR-structures.

Two main approaches are known for the solutions of the equivalence problem: the one of *normal forms* and the one of *differential systems*. The normal form approach works for CR-submanifolds of $\mathbb{C}^N$ (i.e. for embedded CR-structures) and concerns with normalization of the Taylor series of their defining equations at a given point (see e.g. the survey [BER00b] for the references). Once a normal form is established, it solves the local biholomorphic equivalence problem for real-analytic CR-structures but gives neither a solution of the local equivalence problem for smooth CR-structures nor of the global equivalence problem.

In this paper we follow the second approach. Pioneer work here was done by E. CARTAN [CaE32] who gave a solution of the equivalence problem for smooth *Levi-nondegenerate* hypersurfaces in $\mathbb{C}^2$ and later by N. TANAKA [Ta62, Ta67, Ta70] and by S.S. CHERN and J.K. MOSER [CM74] (see also D. BURNS and S. SHNIDER [BS75a, BS75b] and H. JACOBOWITZ [J77]) in $\mathbb{C}^N$, where the problem is reduced to the equivalence problem for the so-called *absolute parallellisms*. The latter means a choice of a frame in each tangent space. The reason that the equivalence problem for absolute parallellisms is well-understood (see e.g. [St83]) is that the condition for a diffeomorphism $f$ to send frames into frames can be reduced to ordinary differential equations for $f$ along smooth real curves. Further reductions to absolute parallellisms have been obtained by R.I. MIZNER [Miz89], T. GARRITY and R.I. MIZNER [GM97], V.V. EZHOV, A.V. ISAEV and G. SCHMALZ [EIS99], A. ČAP and H. SCHICHL [CSi00], G. SCHMALZ and J. SLOVAK [SS00], A. ČAP and G. SCHMALZ [CSm00] for certain Levi-nondegenerate CR-manifolds of higher codimension and by P. EBENFELT [E99] for certain uniformly Levi-degenerate hypersurfaces in $\mathbb{C}^3$.

In the present paper, instead of reducing the problem to the equivalence of absolute parallellisms, we follow a more direct way of reducing the equivalence problem to solving ODE's. This different way allows us to treat more general cases of CR-structures, in particular, all Levi-nondegenerate CR-structures in any codimension. We further treat the case, where the Levi form may change its rank from point to point. To the author's knowledge, no reduction to absolute parallellism is known in this case. A natural example of such a CR-structure is given by a connected oriented hypersurface in $\mathbb{C}^N$ that has both strongly pseudoconvex and strongly pseudoconcave points. Then it is clear from the continuity that there must exist also Levi-degenerate points. It is also clear that such points cannot be avoided by any small perturbation of the hypersurface.

The reduction to ordinary differential equations is obtained by adapting the method of *complete systems* to the equivalence problem. A PDE system is called *complete* if all partial derivatives of the highest order can be expressed through the lower order derivatives. In §1.3 we state Theorem 1.1 establishing a complete system for CR-diffeomorphisms and explain in more detail how it can be used to solve the equivalence problem. Example 1.2 shows that the nondegeneracy conditions in Theorem 1.1 cannot be weakened. In §1.4 an analogous statement gives a solution of the *embedding problem* which is a natural generalization of the equivalence problem. The reader is referred to this paragraph for more detailed discussion of the embedding problem. Then in §2 we establish complete systems in a more general context containing both diffeomorphisms and embeddings as well as some other situations.

The method of complete system has been successfully applied to other problems related to CR-manifolds [Han97, Hay98, Ki99b, E00] and to related problems in differential geometry [Han89, CHY93, Ki99a, CH00]. In fact, the present paper was very much inspired by the work of EBENFELT [E00] (that was, in turn, inspired by the work of C.-K. HAN [Han97]), where complete systems were applied to local jet parametrizations and unique determinations of CR-diffeomorphisms for smooth finitely nondegenerate hypersurfaces in $\mathbb{C}^N$ (see below for the definition). On the other hand, both methods of [Han97] and of [E00] to derive complete systems are very different from ours and yield in their cases systems of higher order than that given by Theorems 1.1 and 2.1 below.

The method to derive complete systems in this paper is based on a development of the *Segre set method* (originally due to M.S. BAOUENDI, P. EBENFELT and L.P. ROTHSCHILD [BER96]). The Segre sets method proved to be a powerful tool for studying the space of (real-analytic) CR-maps between given *real-analytic* CR-submanifolds. The applications of the Segre sets method include a characterization of the finite type [BER96], unique determination and an (explicit algorithmic) holomorphic parametrization of CR-maps by their jets [BER97, Z97, BER99b, Ki99b, L99b], algebraicity of (real-analytic) CR-maps between real-algebraic CR-manifolds [BER96, Z99] and convergence properties of formal maps [BER00a, L99b, BRZ00, Mir00a, Mir00b, BMR01] (see also the references in the quoted papers about related results obtained by other methods).

The construction of the Segre sets can be directly extended neither to abstract non-embeddable (into some $\mathbb{C}^N$) nor to smoothly embeddable CR-structures because it makes use of the complexified defining function: If a submanifold of $\mathbb{C}^N$ is locally defined by a power series $\rho(z,\overline{z}) = 0$, the complexification $\rho(z,\overline{w})$ with $z \neq w$ is needed already for the definition of the *Segre varieties* $Q_w := \{z : \rho(z,\overline{w}) = 0\}$. The Segre sets are then defined inductively by $Q_w^1 := Q_w$, $Q_w^{s+1} := \cup_{z \in Q_w^s} Q_z$ through the Segre varieties. Because of this analytic nature of the Segre sets, the above mentioned applications of the Segre set method cannot be directly extended to the case of smooth CR-structures considered in this paper. In this paper we develop *approximate Segre sets method* that is applied to any smooth CR-submanifold of $\mathbb{C}^N$ and, more generally, to any abstract smooth CR-structure.

We would like to mention that the approximate Segre sets method in this paper yields a complete system in an explicit and algorithmic way. In particular, CR-manifolds of finite smoothness are treated with explicit estimates for the required initial regularity in terms of their "degree of nondegeneracy".

1.2. **Levi-nondegeneracy and higher order nondegeneracies.** We begin with CR-submanifolds of $\mathbb{C}^N$ that serve as basic examples of (abstract) CR-manifolds. A submanifold $M \subset \mathbb{C}^N$ is called CR if its complex tangent space $T_p^c M := T_p M \cap J T_p M$ has constant dimension for $p \in M$, where $J : T\mathbb{C}^N \to T\mathbb{C}^N$ denotes the standard complex structure. A $(1,0)$ (resp. $(0,1)$) vector field on $M$ is a section of the subbundle $T^{1,0}M \subset T^c M \otimes \mathbb{C}$ (resp. $T^{0,1}M \subset T^c M \otimes \mathbb{C}$) consisting of all vectors of the form $X - iJX$ (resp. $X + iJX$) for $X \in T^c M$. If $M$ is of class $\mathcal{C}^\kappa$, its $(1,0)$ and $(0,1)$ vector fields can be considered of class $\mathcal{C}^{\kappa-1}$. The best known invariant for CR-manifolds is the *Levi form* defined to be the (unique) hermitian form

$$\mathcal{L}_p \colon T_p^{0,1}M \times T_p^{0,1}M \to (T_p M/T_p^c M) \otimes \mathbb{C} \tag{1.1}$$

satisfying

$$\mathcal{L}_p(L_1(p), L_2(p)) = \frac{1}{2i}\pi[L_1, \overline{L}_2](p) \tag{1.2}$$

for all $(0,1)$ vector fields $L_1, L_2$, where $\pi \colon TM \otimes \mathbb{C} \to (TM/T^c M) \otimes \mathbb{C}$ is the canonical projection. According to V.K. BELOSHAPKA [Be89], the Levi form $\mathcal{L}_p$ is called *nondegenerate* if (1) $\mathcal{L}_p(L_1, L_2) = 0$ for all $L_2$ implies $L_1 = 0$ and (2) the vectors $\mathcal{L}_p(L_1, L_2)$ for all $L_1$ and $L_2$ span the vector space $(T_p M/T_p^c M) \otimes \mathbb{C}$. In [Be89] BELOSHAPKA showed that, if $M$ is a real quadric of the form $\{\mathsf{Im}\, w = H(z, z)\}$, where $H \colon \mathbb{C}^m \times \mathbb{C}^m \to \mathbb{C}^d$ is hermitian, this condition of nondegeneracy is necessary and sufficient for the local stability group (i.e. for the local biholomorphisms fixing a point) of $M$ to be a (finite-dimensional) Lie group.

For manifolds that are more general than quadrics, the nondegeneracy of the Levi form appears to be only a sufficient condition for the finite-dimensionality of the local stability group. For instance, the tube over the light cone in $\mathbb{C}^N$ given by $(\mathsf{Re}\, z_n)^2 = (\mathsf{Re}\, z_1)^2 + \cdots + (\mathsf{Re}\, z_{n-1})^2$ is everywhere Levi-degenerate but has a finite-dimensional stability group. In order to have optimal conditions, one often has to replace Levi-nondegeneracy by finer *higher order* nondegeneracy conditions (see e.g. the survey [BER00b] for related references). Natural generalizations of the conditions (1) and (2) in the above definition of the nondegenerate Levi form are conditions of finite type and of finite nondegeneracy. Recall that a CR-submanifold $M$ is said to be of *finite type* at a point $p \in M$ (in the sense of J.J. KOHN [Ko72] and T. BLOOM and I. GRAHAM [BG77]) if all $(1,0)$ and $(0,1)$ vector fields on $M$ together with all their higher order commutators span the space $T_p M \otimes \mathbb{C}$. The minimal length $\nu \geq 1$ of the commutators required to span $T_p M \otimes \mathbb{C}$ is called the *type* of $M$ at $p$ and is denoted here by $\nu$. The *type of $M$* (without specifying a point) is the maximal type of $M$ at $p$ for all $p \in M$. A CR-submanifold $M \subset \mathbb{C}^N$ is called *finitely nondegenerate* at $p$ if there exists $l \geq 0$ such that for any $(0,1)$ vector field $L$ on $M$ with $L(p) \neq 0$, there are $(0,1)$ vector fields $L_1, \ldots, L_k$ on $M$, $0 \leq k \leq l$, satisfying

$$[L_1, \ldots, [L_k, \overline{L}] \ldots](p) \notin T_p^c M \otimes \mathbb{C}. \tag{1.3}$$

We call the minimal number $l$ with this property the *degeneracy* of $M$ at $p$ and $M$ to be $l$-*nondegenerate* at $p$. The *degeneracy of $M$* is the maximal degeneracy of $M$ at $p$ for all $p \in M$.

The reader is referred e.g. to [BRZ01a] for further discussion on this condition for real-analytic CR-submanifolds. Both conditions to have a type $\nu$ and a degeneracy $l$ at $p$ involve finite number of differentiations and therefore are meaningful for $\mathcal{C}^\kappa$-smooth CR-submanifolds provided $\nu \leq \kappa$ and $l \leq \kappa - 1$. In contrast to the Levi-nondegeneracy, the condition of finite nondegeneracy and of finite type can be always achieved by small perturbations (see [BRZ01b]).

1.3. **Complete systems for CR-diffeomorphisms.** In this section by saying "smooth" we always mean $\mathcal{C}^\infty$. Recall that an *abstract CR-manifold of CR-dimension $n$ and of CR-codimension $d$* (see e.g. [Bo91, BER99a]) is a real $(2n+d)$-dimensional (smooth) manifold $M$ together with a subbundle $T^c M \subset TM$ of rank $2n$ and a complex structure $J$ viewed as a bundle automorphism $J \colon T^c M \to T^c M$ with $J^2 = -\mathsf{id}$ satisfying the following *integrability condition*: for any two $(0,1)$ vector fields (defined the same way as before), the commutator is also a $(0,1)$ vector field. Here $(1,0)$ and $(0,1)$ vector fields, finite type and finite nondegeneracy are defined by repeating the definitions in the embedded case given above.

Given two manifolds $M$ and $M'$ and an integer $k$, denote by $J^k(M, M')$ the bundle of all $k$-jets of smooth maps between (open pieces of) $M$ and $M'$. If, in addition, $\dim M = \dim M'$, denote by $G^k(M, M') \subset J^k(M, M')$ the open subset of all jets of invertible maps. One of the main results of this paper that gives a solution of the equivalence problem is the existence of complete systems in the jet bundles of CR-diffeomorphisms:

**Theorem 1.1.** *Let $M$ and $M'$ be smooth (abstract) CR-manifolds, both finitely nondegenerate and of finite type. Then there exist a number $r$ and a smooth map $\Phi \colon G^r(M, M') \to J^{r+1}(M, M')$ such that every smooth CR-diffeomorphism $f$ between open pieces of $M$ and $M'$ satisfies the complete differential system*

$$j_x^{r+1} f \equiv \Phi(x, j_x^r f) \tag{1.4}$$

*for all $x$ in the domain of definition of $f$. The number $r$ can be chosen to be $2(d+1)l$, where $d$ is the CR-codimension and $l$ is the degeneracy of $M$.*

*Remark.* Since both conditions of finite type and finite nondegeneracy are invariant under CR-diffeomorphisms, it is sufficient to require them only for $M$ or only for $M'$.

Theorem 1.1 and Corollaries 1.3–1.5 in this paragraph are special cases of the more general statements given by Theorem 2.1 and Corollaries 2.2–2.4 below.

The condition of finite nondegeneracy in Theorem 1.1 cannot be replaced by a weaker known condition, i.e. by essential finiteness (see e.g. [BER99a]) as the following example shows.

*Example* 1.2. The hypersurface

$$M := \{(z, w) \in \mathbb{C}^2 : \mathsf{Im}\, w = |z|^4\}$$

is essentially finite and of finite type but not finitely nondegenerate at 0. We claim that the statement of Theorem 1.4 does not hold for $M$ and $M' := M$. If it held, Corollary 1.5 below would imply that the dimension of the local stability group $\mathsf{Aut}(M, q)$ is upper-semicontinuous with respect to $q \in M$. (Recall that $\mathsf{Aut}(M, q)$ consists of all germs at $q$ of biholomorphic maps of $\mathbb{C}^2$ fixing $q$ and sending $M$ into itself.) On the other hand, $M$ is locally biholomorphically equivalent to the quadric $Q := \{\operatorname{Im} w = |z|^2\}$ at any point $(z, w)$ with $z \neq 0$ via the map $(z, w) \mapsto (z^2, w)$. Hence the group $\operatorname{Aut}(M, q)$ can be computed directly and it is easy to see that its dimension is not upper-semicontinuous which is the contradiction proving the claim.

In the case where $M$ and $M'$ are hypersurfaces in $\mathbb{C}^N$, a similar conclusion has been recently obtained by EBENFELT [E00] with the order $r = l^3 + l^2 + l + 1$ which is larger than $r = 4l$ given by Theorem 1.1 in the hypersurface case ($d = 1$).

Given the complete system (1.4), a solution of the equivalence problem for $M$ and $M'$ can be obtained as follows. For each initial value of the jet $\lambda = j_x^r f$ at a given point $x \in M$, the whole map $f$ can be uniquely reconstructed by solving ODE's along real curves in $M$. Then $M$ and $M'$ are equivalent if and only if there exists an initial value $\lambda \in G^r(M, M')$ such that the corresponding map $f$ with $j_x^r f = \lambda$ extends to a CR-equivalence between $M$ and $M'$.

Recall that, given a fiber bundle $\pi \colon F \to B$, a *connection* on $\pi$ is a subbundle $H$ of the tangent bundle $TF$ such that the restriction of $d\pi$ on each $H_x \subset T_x F$, $x \in F$, is an isomorphism onto $T_{\pi(x)} B$. Given a connection $H$, a section $s \colon U \to F$ of $\pi$ (with $U \subset B$ open) is called *horisontal* if $ds(T_y B) = H_{s(y)}$ for every $y \in U$. Under the above nondegeneracy assumptions we obtain a differential-geometric solution of the equivalence problem:

**Corollary 1.3.** *Suppose that the conclusion of Theorem 1.1 holds. Then there exists a smooth connection on the fiber bundle $\pi \colon G^r(M, M') \to M$ such that a smooth section $s$ of $\pi$ is horisontal if and only if $s(x) \equiv j_x^k f$ for some CR-diffeomorphism $f$ between open pieces of $M$ and $M'$.*

Thus the equivalence problem is reduced to finding horisontal sections of the given connection. Another application of Theorem 1.1 is the following unique determination result:

**Corollary 1.4.** *Suppose that the conclusion of Theorem 1.1 holds. Then CR-diffeomorphisms between connected open pieces of $M$ and $M'$ are uniquely determined by their $r$-jets at any point, i.e. if $f, g \colon U \subset M \to M'$ are CR-diffeomorphisms, $U$ is connected, $p \in U$ and $j_p^r f = j_p^r g$, then $f \equiv g$.*

Furthermore, we obtain a local jet-parametrization, meaning that the CR-diffeomorphisms are not only uniquely determined by their $r$-jets but also depend on them in a smooth fashion. We denote by $G_{p,p'}^r(M, M')$ the fiber in $G^r(M, M')$ over $(p, p') \in M \times M'$, i.e. the $r$-jets of maps with fixed source $p$ and target $p'$.

**Corollary 1.5.** *Suppose that the conclusion of Theorem 1.1 holds. Then for any $p \in M$ and $p' \in M'$, there exist an open neighborhood $\Omega$ of $\{p\} \times G_{p,p'}^r(M, M')$ in $M \times G^r(M, M')$ and a smooth map $\Psi \colon \Omega \to M'$ such that every smooth CR-map $f$ from a neighborhood of a point $q \in M$ into $M'$ with $(q, j_q^r f) \in \Omega$ satisfies the identity*

$$f(x) \equiv \Psi(x, q, j_q^r f) \tag{1.5}$$

*for all $x \in M$ sufficiently close to $q$.*

We remark that Corollary 1.5 gives a stronger statement than in [BER97, Z97, BER99b, Ki99b, L99b] in the cases considered there. The main difference is that the parametrization in Corollary 1.5 is obtained simultaneously for all points $q$ near a given one and for all local CR-maps in a neighborhood of $q$. The idea of the proof of Corollary 1.5 is to integrate the system (1.4) along curves in $M$ chosen in a proper way. Vice versa, given a parametrization (1.5), the statement of Theorem 1.1 can be obtained by differentiating (1.5). Here the mentioned stronger version of a parametrization with the varying point $q$ is essential.

Using recent results of [BRWZ01] we obtain, as an application of Corollary 1.5, a Lie group structure on the group $\mathsf{Aut}(M)$ of all (global) smooth CR-automorphisms of $M$. We view $\mathsf{Aut}(M)$ as a closed subgroup in the topological group $\mathsf{Diff}(M)$ of all smooth diffeomorphisms of $M$. Here a sequence in $\mathsf{Diff}(M)$ converges to an element there if it converges uniformly on compacta with all derivatives.

**Corollary 1.6.** *Let $M$ be a smooth abstract CR-manifold that is finitely nondegenerate and of finite type. Then $\mathsf{Aut}(M)$ is a Lie transformation group. If $M$ is in addition real-analytic, then the action of $\mathsf{Aut}(M)$ on $M$ is also real-analytic.*

In contrast to the global automorphism group $\mathsf{Aut}(M)$, the behaviour of the local stability group $\mathsf{Aut}(M, p)$ (that consists of germs at $p \in M$ of all local CR-automorphisms of $M$ fixing $p$) is very different for smooth and for real-analytic CR-manifolds. In the real-analytic case, under the same assumptions as in Corollary 1.6, it was shown in [Z97] that $\mathsf{Aut}(M, p)$ is a Lie group whose local action on $M$ is real-analytic. In the smooth case $\mathsf{Aut}(M, p)$ has in general no Lie group structure (see [Z01]).

Finally we give a deformation version of Theorem 1.1 that follows directly from Theorem 4.1 below.

**Theorem 1.7.** *Under the assumptions of Theorem 1.1 suppose that $M = M(t)$ and $M' = M'(t)$ depend smoothly on a parameter $t \in \mathbb{R}^m$ near the origin and fix a smooth identification $G^r(M(t), M'(t)) \cong G^r(M(0), M'(0))$. Then, for every $p \in M(0)$, $p' \in M'(0)$ and $\Lambda \in G^r_{p, p'}(M(0), M'(0))$, there exist a neighborhood $U$ of $(p, p', \Lambda)$ in $G^r(M(0), M'(0))$, and a smooth family of maps $\Phi_t \colon U \to G^{r+1}(M(0), M'(0))$ such that, for $t \in \mathbb{R}^m$ near the origin, every smooth CR-diffeomorphism $f$ from a neighborhood of a point $q \in M(t)$ into $M'(t)$ with $(q, j^r_q f) \in U$ satisfies $j^{r+1}_x f \equiv \Phi_t(x, j^r_x f)$. If $M$ and $M'$ are real-analytic, the family $\Phi_t$ can be also chosen to be real-analytic.*

1.4. **Complete systems for CR-embeddings.** A fundamental problem which is a natural generalization of the equivalence problem is the *embedding problem*:

*Given two manifolds with CR-structures, when does there exist a CR-embedding of one manifold into the other?*

The CR-embedding problem is related to the existence of proper holomorphic embeddings between smoothly bounded domains of different dimension the same way the CR-equivalence problem is related to the biholomorphic equivalence problem. Here the existence and regularity properties of embeddings are much less understood (see the survey [Fo93] for related results). Complete systems in this case were constructed for real-analytic hypersurfaces due to [Han97, Hay98, Ki99b].

In this paragraph we show how to extend the results from the previous paragraph to treat CR-embeddings between smooth (abstract) CR-manifolds (of any codimension). We begin by adapting the finite nondegeneracy condition. An embedding of an (abstract) CR-manifold $M$ into an (abstract) CR-manifold $M'$ is called *finitely nondegenerate* at a point $p \in M$ if there exists $l \geq 0$ such that for any $(0,1)$ vector field $L$ on $M'$ with $L(p) \neq 0$, there are $(0,1)$ vector fields $L_1, \ldots, L_k$ on $M'$, $0 \leq k \leq l$, that restrict to $(0,1)$ vector fields on $M$ satisfying

$$[L_1, \ldots, [L_k, \overline{L}] \ldots ](p) \notin T_p^c M' \otimes \mathbb{C}. \tag{1.6}$$

If $l$ is the minimal number with this property, we call the embedding $l$-*nondegenerate*. It is easy to check that if an embedding of $M$ (into some $M'$) is finitely nondegenerate, then $M$ is itself finitely nondegenerate. The converse is not true: the simplest example is the standard linear embedding of the sphere $M = S^3 \subset \mathbb{C}^2$ into the sphere $M' = S^5 \subset \mathbb{C}^3$.

In Theorem 1.1 we had a complete system defined for all invertible jets of smooth maps between $M$ and $M'$ whereas the nondegeneracy conditions were put on $M$ and $M'$. In the case of embeddings we have the same finite type condition on $M$ but the condition of finite nondegeneracy on $M'$ is replaced by finite nondegeneracy of an embedding. In fact, the last condition involves only finitely many derivatives and hence can be considered as a condition on the jet of an embedding. We write here $D^{r,l}(M, M')$ for the (open) subset of all $r$-jets of smooth embeddings that are $\widetilde{l}$-nondegenerate for some $\widetilde{l} \leq l$. (The reader is referred to the next paragraph for a more precise definition and discussion of this notion). We then have the following version of Theorem 1.1 for CR-embeddings:

**Theorem 1.8.** *Let $M$ and $M'$ be smooth (abstract) CR-manifolds, where $M$ is of finite type. Then, for any $l \geq 0$ there exist a number $r > 0$ and a smooth map $\Phi \colon D^{r,l}(M, M') \to J^{r+1}(M, M')$ such that every smooth $l$-nondegenerate CR-embedding $f$ of an open piece of $M$ into $M'$ satisfies the complete differential system*

$$j_x^{r+1} f \equiv \Phi(x, j_x^r f) \tag{1.7}$$

*for all $x$ in the domain of definition of $f$. The number $r$ can be chosen to be $2(d+1)l$, where $d$ is the CR-codimension of $M$.*

Theorem 1.8 is obtained as a special case of the more general statement given by Theorem 2.1 below. As a consequence, we have a solution of the embedding problem for finitely nondegenerate embeddings and the results analogous to Corollaries 1.3–1.5 for them.

A prototype of the above condition of $l$-nondegeneracy (in connection with a different problem) appeared in a paper of J.A. CIMA, S.G. KRANTZ and T.J. SUFFRIDGE [CKS84] and was later improved by J.J. FARAN [Fa90]. FARAN gave his condition in different terms but it can be shown to be equivalent to that of 2-nondegeneracy in his case, where $M \subset \mathbb{C}^N$ and $M' \subset \mathbb{C}^{N'}$ are strongly pseudoconvex (or Levi-nondegenerate) real-analytic hypersurfaces. Later C.-K. HAN [Han97] gave a condition equivalent to that of $l$-nondegeneracy in this case. In all these cases a certain partial normalization of the defining equation of $M'$ has been used (that is known to exist due to the CHERN-MOSER normal form [CM74]). Note that such a normalization is not possible in general if $M'$ is not a Levi-nondegenerate hypersurface. More recently the condition of $l$-nondegeneracy was extended by the first author [Ki99b] to general real-analytic hypersurfaces and by B. LAMEL [L99a, L99b, L00] to smooth embedded CR-manifolds of any codimension.

Remarkably, the condition of finite nondegeneracy for CR-embeddings that may look restrictive at the first glance, holds sometimes automatically for all embeddings of manifolds of certain classes. We give two situations where this phenomenon occurs. The first situation was considered by S.M. WEBSTER [W79] in the case of a smooth CR-manifold $M$ of hypersurface type embedded into the unit sphere $M' = S^{2n+1} \subset \mathbb{C}^{n+1}$ such that $M$ has the real codimension 2 in $M'$. In this setting WEBSTER established a fundamental relation between the pseudo-conformal curvature tensor $S_M$ of $M$ defined by CHERN and MOSER [CM74] and the second fundamental form of the embedding of $M$ into $M'$. In fact, it follows from his formula (2.14) that if $S_M$ does not vanish at a point $p \in M$, then the restriction of the fundamental form to $T_p^c M$ also does not vanish for any embedding of $M$ into $M'$. But the last condition means exactly 2-nondegeneracy of the embedding as was shown by FARAN [Fa90]. Hence we obtain the following application of Theorem 1.8:

**Corollary 1.9.** *Let $M' = S^{2n+1} \subset \mathbb{C}^{n+1}$ be the unit sphere and $M$ be any smooth strongly pseudoconvex CR-manifold of hypersurface type of dimension $2n-1$ whose tensor $S_M$ does not vanish at any point. Then the conclusion of Theorem 1.8 holds for all CR-embeddings of an open piece of $M$ into $M'$ with $r = 8$.*

Recall that $S_M$ always vanish when $M$ is 3-dimensional. On the other hand, in any higher dimension, if $(M, p)$ is not locally equivalent to the sphere for some $p \in M$, then $S_M$ does not identically vanish in a neighborhood of $p$.

We conclude by giving another general situation where all embeddings are automatically finitely nondegenerate. This time a condition is put only on the target manifold $M'$. We call a smooth CR-manifold $M'$ *finitely nondegenerate in dimension $n$* at a point $p \in M'$ if for any linearly independent $(0,1)$ vector fields $L_1, \ldots, L_n$ on $M'$ near $p$ and another $(0,1)$ vector field $L$ with $L(p) \neq 0$, there exist a number $k \geq 0$ and indices $j_1, \ldots, j_k \in \{1, \ldots, n\}$ such that

$$[L_{j_1}, \ldots, [L_{j_k}, \overline{L}] \ldots](p) \notin T_p^c M' \otimes \mathbb{C}. \tag{1.8}$$

The notion of $l$-nondegeneracy in dimension $n$ is defined in the obvious way as above. It follows directly from the definition that if $M'$ is finitely nondegenerate in dimension $n$ at a point $p'$ and if $M$ has CR-dimension at least $n$ then any embedding of $M$ into $M'$ through $p'$ is finitely nondegenerate at this point. The following example shows that nondegenerate manifolds in the above sense exist already in the lowest dimension.

*Example* 1.10. We claim that the hypersurface $M' \subset \mathbb{C}^3$ given by

$$M' := \{(z_1, z_2, w) : \operatorname{Im} w = |z_1|^2 + |z_2|^2 + \operatorname{Im}\left(z_1^2 \overline{z}_1 + (z_1 + z_2)^3(\overline{z}_1 + \overline{z}_2)\right)\}$$

is 3-nondegenerate in dimension 1 at $p = 0$. Thus any embedding of any hypersurface $M \subset \mathbb{C}^2$ into $M'$ through 0 is automatically $l$-nondegenerate at 0 for some $l \leq 3$. Indeed, given any $(0,1)$ vector fields $L_1$ and $L$ with $L_1(0) \neq 0$, $L(0) \neq 0$, one has $[L_1, \overline{L}](p) \in T_p^c M' \otimes \mathbb{C}$ if and only if $L_1(0)$ and $L(0)$ are orthogonal with respect to the Levi form $|z_1|^2 + |z_2|^2$. If, in addition $[L_1, [L_1, \overline{L}]](p) \in T_p^c M' \otimes \mathbb{C}$, then either $L_1(0) = (a, 0)$ and $L(0) = (0, b)$ or $L(0) = (a, 0)$ and $L_1(0) = (0, b)$ for some $a, b \neq 0$. But in this case one always has $[L_1, [L_1, [L_1, \overline{L}]]](p) \notin T_p^c M' \otimes \mathbb{C}$, which proves the claim.

We now have the following consequence of Theorem 1.8:

**Theorem 1.11.** *Let $M$ and $M'$ be smooth (abstract) CR-manifolds, $M$ being of finite type and CR-dimension $n$ and CR-codimension $d$ and $M'$ being $l$-nondegenerate in dimension $n$. Then the conclusion of Theorem 1.8 holds for all CR-embeddings of an open piece of $M$ into $M'$ with $r = 2(d+1)l$.*

## 2. STATEMENTS OF MAIN RESULTS

We now allow $M$ and $M'$ to be abstract CR-manifolds of possibly different CR-dimensions and CR-codimensions and consider CR-maps of finite smoothness that are not necessarily embeddings. If $M$ is of class $\mathcal{C}^\kappa$, the condition on $M$ to have the (finite) type $\nu$ at a point $p \in M$ makes sense for $\nu \leq \kappa$. We further extend the finite nondegeneracy condition from §1.4 to arbitrary CR-maps. As before, this condition will require a restricted number of differentiations and hence will make sense in the case where the map $f$ and the manifolds $M$ and $M'$ have finite degree of smoothness. For an $n$-tuple of integers $\alpha = (\alpha_1, \ldots, \alpha_n)$, we write $|\alpha| := \alpha_1 + \cdots + \alpha_n$ and $L^\alpha := L_1^{\alpha_1} \ldots L_n^{\alpha_n}$, where $L_1, \ldots, L_n$ is a basis of $(0,1)$ vector fields near $p$. Suppose for a moment that $M'$ is embedded in $\mathbb{C}^{N'}$ with a local defining function $\rho' = (\rho'^1, \ldots, \rho'^{d'})$ near a point $p' \in M'$. A $\mathcal{C}^r$-smooth CR-map $f$ between $M$ and $M'$ with $f(p) = p'$ is said to be *$l$-nondegenerate* at $p$ if

$$\mathrm{span}_{\mathbb{C}}\{L^\alpha \rho'^j_{Z'}(f(Z), \overline{f(Z)})(p) : 1 \leq j \leq d', 0 \leq |\alpha| \leq l\} = \mathbb{C}^{N'}, \tag{2.1}$$

where $\rho'^j_{Z'} = (\partial \rho'^j / \partial Z_1, \ldots, \partial \rho'^j / \partial Z_{N'})$ denotes the complex gradient, and if $l \geq 0$ is the minimal number for which (2.1) holds. The definition makes sense for $l \leq \min(r, \kappa - 1)$. In this form it is due to [Ki99b, L99a, L99b]. The arguments in [L99a, L99b] show that this definition does not depend on the choice of the basis $L_1, \ldots, L_n$, the defining function $\rho'$ and the coordinates $Z'$ (where the gradient $\rho'^j_{Z'}$ is computed). In particular, the given definition does not depend on the embedding of $M$ into $\mathbb{C}^{N'}$.

An inspection of the above definition shows that it does not really use the whole CR-map $f$ but rather its $l$-jet at $p$. Moreover, condition (2.1) can be written for any map $f \colon M \to M'$, not necessarily CR. It will be still independent of the choice of $L_1, \ldots, L_n$ and of $\rho'$ but will in general depend on the coordinates $Z'$. In order to eliminate this dependence, we restrict ourselves in this definition to the $r$-jets that we call *CR $r$-jets* and that are defined to be those $r$-jets at a point $p \in M$ of $\mathcal{C}^r$-smooth maps $f \colon M \to M'$ that satisfy $Lf(p) = o(|x-p|^{r-1})$ for any $(0,1)$ vector field $L$ on $M$. Any jet of a CR-map has this property but, in general, there might be more CR $r$-jets at $p$ than CR-maps defined near $p$. We now call a CR $r$-jet $(p, \Lambda) \in J^r(M, M')$ *$l$-nondegenerate* $(l \leq r)$ if (2.1) holds for some (and hence for any) representative $f$ of $(p, \Lambda)$. Again, the arguments of [L99a, L99b] show that this definition is independent of the choice of $L_1, \ldots, L_n$, $\rho'$ and $Z'$.

We now extend the notion of $l$-nondegenerate maps to the case, where $M'$ is not necessarily embedded in the same spirit as it was done in §1.2 for $l$-nondegenerate CR-manifolds and in §1.4 for CR-embeddings. Since we wish to allow CR-maps that are not embeddings we cannot in general push $(0,1)$ vector fields from $M$ to $M'$. The trick to avoid this difficulty is to push these vector fields to $M \times M'$ via the graph of $f$. In the product $M \times M'$ we identify a vector field $L'$ on $M'$ with the "vertical" vector field $\widetilde{L'}(x, x') := (0, L'(x'))$ and a vector field $L$ on $M$ with the "pushed" vector field $\widetilde{L}(x, f(x)) := (L(x), f_* L(x))$ along the graph of $f$. We now call a $\mathcal{C}^r$-smooth CR-map $f$ as above $l$-*nondegenerate* at $p$ if, for any $(0, 1)$ vector field $L'$ on $M'$ with $L'(p') \neq 0$, there are $(0, 1)$ vector fields $L_1, \ldots, L_k$ on $M$, $0 \leq k \leq l$, and their extensions $\widetilde{L}_1, \ldots, \widetilde{L}_k$ from the graph of $f$ to $(0, 1)$ vector fields on $M \times M'$ such that

$$[\widetilde{L}_1, \ldots, [\widetilde{L}_k, \overline{L'}] \ldots](p, p') \notin (T_p M \oplus T_{p'}^c M') \otimes \mathbb{C} \tag{2.2}$$

and $l \geq 0$ is the minimal number with this property. We call $l$ the *degeneracy* of $f$ at $p$. If $f$ is defined on an open subset $U \subset M$, we call it $l$-nondegenerate if $l$ is the maximal degeneracy of $f$ at $p$ for all $p \in U$. Characterizations of $l$-nondegeneracy as in [BRZ01a, §3] can be also extended to this more general case in an analogous way. In particular, if $M'$ is embedded, the above two definitions are equivalent. On the other hand, if $M'$ is merely abstract, it can be "approximately embedded" of sufficiently high order in the sense of Proposition 3.1 below and then $l$-nondegeneracy of CR-maps into $M'$ can be defined through any such embedding. Finally, the notions of CR $r$-jets and their $l$-nondegeneracy introduced above can be also directly extended to abstract CR-manifolds in the obvious way.

As the main result of this section, we obtain a complete differential system on the jet bundle $J^r(M, M')$ that is satisfied by all $l$-nondegenerate CR-maps between $M$ and $M'$. We express it as a section of the fiber bundle $J^{r+1}(M, M') \to J^r(M, M')$. A natural replacement of the open subset $G^r(M, M') \subset J^r(M, M')$ of all invertible jets (considered in Theorem 1.1) is the invariantly defined open subset of all those $r$-jets in $J^r(M, M')$ that are either not CR or CR and $\widetilde{l}$-nondegenerate for some $0 \leq \widetilde{l} \leq l$. We denote this subset by $D^{r,l}(M, M') \subset J^r(M, M')$.

**Theorem 2.1.** *Let $M$ and $M'$ be $\mathcal{C}^\kappa$-smooth abstract CR-manifolds and denote by $d$ the CR-codimension of $M$. Suppose that $M$ is of the type $2 \leq \nu \leq \kappa$, fix a number $l \geq 0$ and set $r := 2(d+1)l$ and*

$$k := 4(d^2 + d)\nu l + 4(d^2 - 1)l + 2d\nu - 2d + 1. \tag{2.3}$$

*Then, if $\kappa \geq k + 1$, there exists a $\mathcal{C}^{\kappa - k - 1}$-smooth section $\Phi \colon D^{r,l}(M, M') \to J^{r+1}(M, M')$ such that, for every open subset $U \subset M$ and every $0 \leq \widetilde{l} \leq l$, every $\widetilde{l}$-nondegenerate $\mathcal{C}^k$-smooth CR-map $f \colon U \to M'$ satisfies*

$$j_x^{r+1} f \equiv \Phi(x, j_x^r f), \quad x \in U. \tag{2.4}$$

In the special case when $M$ is of hypersurface type (i.e. $d = 1$), (2.3) yields $k = 8\nu l + 2\nu - 1$. In general, a more refined estimate can be obtain by replacing $\nu$ with the sum of all Hörmander numbers of $M$ at $p$ due to Lemma 3.4 and Theorem 4.1 below.

*Remark.* It follows from the proof of Theorem 2.1 that the same conclusion holds also for every $x \in M$ and every $l$-nondegenerate $\mathcal{C}^k$-smooth map $f$ from a neighborhood of $x$ in $M$ into $M'$ that is only CR up to order $k$ at $x$ rather than "precisely" CR.

We use Theorem 2.1 to obtain a smooth local parametrization of CR-maps by their $r$-jets as follows:

**Corollary 2.2.** *Suppose that the conclusion of Theorem 2.1 holds. Then, for every $(x_0, \lambda_0) \in D^{r,l}(M, M')$, there exists a $\mathcal{C}^{\kappa - k - 1}$-smooth map $\Theta$ from a neighborhood $\Omega$ of $(x_0, x_0, \lambda_0)$ in $M \times$* $J^k(M, M')$ into $M'$ such that for every $p \in M$, every $\mathcal{C}^k$-smooth CR-map $f$ from a neighborhood $U(p) \subset M$ into $M'$ with $(p, p, j_p^r f) \in \Omega$ satisfies the identity

$$f(x) \equiv \Theta(x, p, j_p^r f)$$

for all $x \in M$ in a neighborhood $\widetilde{U}(p) \subset U(p)$ that depends on $U(p)$ and $\lambda_0$ but not on $f$.

An immediate consequence of Corollary 2.2 is the following regularity property for CR-maps:

**Corollary 2.3.** *Suppose that the conclusion of Theorem 2.1 holds. Then any l-nondegenerate $\mathcal{C}^k$-smooth CR-map between open pieces of $M$ and $M'$ is necessarily of class $\mathcal{C}^{\kappa-k-1}$. In particular, if $M$ and $M'$ are $\mathcal{C}^\infty$, then $f$ is also $\mathcal{C}^\infty$.*

For $\kappa = \infty$ and, when $M$ and $M'$ are embeddable, the statement of Corollary 2.3 was recently obtained by LAMEL [L00]. As a consequence of Theorem 2.1, we obtain a differential-geometric description of all possible locally defined nondegenerate CR-maps via horisontal sections in the jet bundle $J^r(M, M')$:

**Corollary 2.4.** *Suppose that the conclusion of Theorem 2.1 holds. Then there exists a $\mathcal{C}^{\kappa-k-1}$-smooth connection on the fibration $D^{r,l}(M, M') \to M$ such that a $\mathcal{C}^{k-r}$-smooth section $s \colon U \subset M \to D^{r,l}(M, M')$ is horisontal if and only if $s(x) \equiv j_x^r f$ for a $\mathcal{C}^k$-smooth CR-map $f \colon U \to M'$.*

## 3. Approximate Segre sets and iterated complexifications

**3.1. Approximation of abstract CR-manifolds by embedded ones.** Our first main tool is the approximation of abstract CR-manifolds by embedded real-analytic ones. It is well-known (see e.g. [Bo91]) that there exists a smooth abstract CR-manifold that is not locally embeddable in any $\mathbb{C}^N$. On the other hand, even if a CR-manifold $M$ is smoothly embedded into $\mathbb{C}^N$, it can be shown that, in general, it may not be embedded as a real-analytic CR-submanifold [Z01]. Here we have a different point of view: we look for embeddings that are CR only up to some order. The following proposition shows that such embeddings always exist.

**Proposition 3.1.** *Let $M$ be a $\mathcal{C}^\kappa$-smooth abstract CR-manifold $(2 \le \kappa \le \infty)$ of CR-dimension $n$ and CR-codimension $d$ and let $k \le \kappa$ be any integer. Then, for any $p \in M$, there exists a neighborhood $U$ of $p$ in $M$ and a $\mathcal{C}^{\kappa-k}$-smooth map $\varphi \colon U \times U \to \mathbb{C}^{n+d}$ such that the following hold:*

(i) *for every $x \in U$, $\varphi_x := \varphi(x, \cdot)$ is an embedding of $U$ into $\mathbb{C}^{n+d}$ whose image is a generic real-analytic CR-submanifold of codimension $d$;*

(ii) *the map $\varphi_x$ is CR of order $k$ at $x$, i.e. $L\varphi_x(y) = o(|x - y|^{k-1})$ as $y \to x$ for any $(0, 1)$ vector field $L$ on $M$ defined in some neighborhood of $x$.*

For the proof of Proposition 3.1 we need the following approximate version of the holomorphic Frobenius theorem:

**Lemma 3.2.** *Let $k \ge 0$ be any integer and $L_1, \ldots, L_n$ be linearly independent holomorphic vector fields in $\mathbb{C}^m$ defined in a neighborhood of $0$ such that*

$$[L_i, L_j](Z) = \sum c_{i,j}^s(Z) L_s(Z) + o(|Z|^k), \quad Z \to 0, \quad 1 \le i, j \le n, \tag{3.1}$$

*for some holomorphic functions $c_{i,j}^s(Z)$. Then there exist holomorphic vector fields $\widetilde{L}_1, \ldots, \widetilde{L}_n$ spanning the same subspace as $L_1, \ldots, L_n$ at every point and a local holomorphic change of coordinates $\widetilde{Z}$ such that*

$$\widetilde{L}_i = \frac{\partial}{\partial \widetilde{Z}_i} + o(|\widetilde{Z}|^{k+1}), \quad \widetilde{Z} \to 0, \quad 1 \le i \le n.$$

*Moreover, if $L_1, \ldots, L_n$ depend smoothly on real parameters, the vector fields $\widetilde{L}_1, \ldots, \widetilde{L}_n$ and the coordinates $\widetilde{Z}$ can be chosen to depend smoothly on the same parameters with the same smoothness.*

*Proof.* Without loss of generality, $L_i(0) = \frac{\partial}{\partial Z_i}$ for $1 \le i \le n$. We write $Z = (z, w) \in \mathbb{C}^n \times \mathbb{C}^{m-n}$ for the local coordinates. Then we can further assume that $L_i(z, w) = \frac{\partial}{\partial z_i} + \sum a_i^j(z, w) \frac{\partial}{\partial w_j}$ for some holomorphic functions $a_i^j(z, w)$. For this choice of vector fields, (3.1) implies

$$[L_i, L_j](Z) = o(|Z|^k), \quad Z \to 0, \quad 1 \le i, j \le n. \tag{3.2}$$

It is easy to see that there exists a holomorphic map $\chi$ from a neighborhood of 0 in $\mathbb{C}^m$ into $\mathbb{C}^{m-n}$ such that

$$(\partial_z^\alpha \chi^i)(0, w) = (L^\alpha w^i)(0, w), \quad 0 \le |\alpha| \le k+2, \quad 1 \le i \le m-n, \tag{3.3}$$

where $\alpha \in \mathbb{Z}_+^n$ is any multiindex. Furthermore, the identity (3.2) allows us to permute the order of applying the vector fields on the right-hand side of (3.3) with a controlled error term:

$$(\partial_{z_{i_1}} \ldots \partial_{z_{i_s}} \chi^i)(0, w) = (L_{i_1} \ldots L_{i_s} w^i)(0, w) + o(|w|^{k+2-s}), \tag{3.4}$$

where $i_1, \ldots, i_s \in \{1, \ldots, n\}$, $0 \le s \le k+2$, $1 \le i \le m-n$. Then it follows that the change of coordinates given by $\widetilde{Z}(z, w) = (z, \chi(z, w))$ satisfies the required conclusion. The last statement about the dependence on parameters follows directly from the proof. $\square$

*Proof of Proposition 3.1.* The proof can be obtained by following the arguments in the proof of Theorem 2.1.11 in [BER99a] and using approximate complexifications of the smooth vector fields and Lemma 3.2 instead of Theorem 2.1.12 in [BER99a]. The details are left to the reader. $\square$

### 3.2. Approximate complexifications of abstract CR-manifolds.

In the sequel we use the semicolon ";" to separate the variables, where the functions (or maps) are holomorphic, from the other variables, where they are merely continuously differentiable of certain order.

Let $M$ be an abstract $\mathcal{C}^\kappa$-smooth CR-manifold of CR-dimension $n$ and CR-codimension $d$. We set $N := n + d$. Let $U$ and $\varphi$ be given by Proposition 3.1. We call $M(x) := \varphi_x(U) \subset \mathbb{C}^N$ an *approximate embedding of $M$ of order $k$ at $x$.* By a simple translation argument we may assume that $\varphi_x(x) = 0$ for all $x \in U$. Denote by $\rho(Z, \overline{Z}; x)$ a family of local ($\mathbb{C}^d$-valued) defining functions of $M(x)$ in a neighborhood of 0 that are holomorphic in the first two arguments and $\mathbb{C}^{\kappa-k}$-smooth in all the arguments. As usual, we assume $d\rho^1(\cdot, \cdot; x) \wedge \cdots \wedge d\rho^d(\cdot, \cdot; x) \ne 0$ on the domain of definition for each fixed $x$. Finally, we can choose a neighborhood, where $\rho$ is defined, to be independent of $k$.

We now define an *approximate complexification of $M$ at $x$ of order $k$* by

$$\mathcal{M}(x) := \{(Z, \zeta) : \rho(Z, \zeta; x) = 0\}.$$

Then $\mathcal{M}(x)$ is a complex submanifold of $\mathbb{C}^N \times \mathbb{C}^N$ through 0 of complex codimension $d$. We call the set

$$\mathcal{M} := \{(Z, \zeta; x) : \rho(Z, \zeta; x) = 0\}$$

a *family of approximate complexifications of $M$ of order $k$*. This family clearly depends on the choice of the approximating family $M(x)$.

We further define, for every $s \geq 1$, an *iterated approximate complexification of $M$ at $x$ of order $k$* by

$$\mathcal{M}^s(x) := \{\xi = (\xi^0, \cdots, \xi^s) \in \mathbb{C}^{(s+1)N} : \rho_j(\xi^{j-1}, \xi^j; x) = 0 \text{ for all } j = 1, \ldots, s\},$$

where

$$\rho_j(\xi^{j-1}, \xi^j; x) := \begin{cases} \overline{\rho}(\xi^{j-1}, \xi^j; x) & \text{if } j \text{ is even} \\ \rho(\xi^{j-1}, \xi^j; x) & \text{if } j \text{ is odd}. \end{cases}$$

We also consider the corresponding *family of iterated approximate complexifications of $M$ of order $k$* defined by

$$\mathcal{M}^s := \{(\xi; x) \in \mathbb{C}^{(s+1)N} \times M : \xi \in \mathcal{M}^s(x)\},$$

(We write the index of $\rho$ below to avoid the confusion with the component indices that we wrote above). It follows that, for every $x \in M$ (close to $p$), $\mathcal{M}^s(x) \subset \mathbb{C}^{(s+1)N}$ is a complex submanifold of dimension $sn + N$ through the origin.

### 3.3. Finite type for smooth abstract CR-manifolds.

We now suppose that $M$ has the finite type $\nu$ at $p$ with $1 \leq \nu \leq k \leq \kappa$. Then, for $x$ close to $p$, the approximation $M(x)$ of order $k$ has also the finite type $\nu$ at $0 \in M(x)$.

By a linear change of coordinates in $\mathbb{C}^N$ and shrinking the domain of definition of $\rho$, if necessary, we may assume that the projection

$$\mathcal{M}(x) \subset \mathbb{C}^n \times \mathbb{C}^d \times \mathbb{C}^N \to \mathbb{C}^n \times \{0\} \times \mathbb{C}^N$$

is biholomorphic onto its image. Using this projection we introduce the coordinates on

$$\mathcal{M}^s(x) \subset \mathbb{C}^{(s+1)N} = \mathbb{C}^N \times \cdots \times \mathbb{C}^N$$

as follows. In each factor $\mathbb{C}^N$ except the last one we take the first $n$ coordinates (among the $N$ available ones) and in the last factor $\mathbb{C}^N$ we take all the $N$ coordinates. We keep the notation $\xi^s \in \mathbb{C}^N$ for the last $N$ coordinates and denote the other coordinates by

$$t = (t^0, \ldots, t^{s-1}) \in \mathbb{C}^n \times \cdots \times \mathbb{C}^n = \mathbb{C}^{ns}.$$

Denote by $(v^0, \ldots, v^s) \in \mathcal{M}^s(x) \subset \mathbb{C}^{(s+1)N}$ the point whose coordinates are $(t, 0)$. Then each component $v^j(t; x) \in \mathbb{C}^N$ is a function of class $\mathcal{C}^{\kappa-k}$ in its variables and holomorphic in $t$ for $x$ fixed.

We now observe that, if the coordinates in $\mathbb{C}^N = \mathbb{C}^n \times \mathbb{C}^d$ are chosen to be *normal* for the real-analytic submanifold $M(p) \subset \mathbb{C}^N$ at $p$ (see e.g. [BER99a], §4.2), then the parametrization $v^s(\cdot, 0; x)$ coincides with the $s$th *Segre map* (or its conjugate) in the sense of [BER99b]. Hence by Theorem 3.1.9 in [BER99b], the generic rank of $v^s(\cdot; x)$ equals $N$ for $s = d + 1$.

By Lemma 4.1.3 in [BER99b], we have

$$v^{2s}(0, t^1, \ldots, t^{s-1}, t^s, t^{s-1}, \ldots, t^1; x) \equiv 0$$

and the generic rank of the partial derivative matrix

$$\frac{\partial v^{2s}}{\partial(t^0, t^{s+1}, t^{s+2}, \ldots, t^{2s-1})}$$

along the linear subspace

$$\{(0, t^1, \ldots, t^{s-1}, t^s, t^{s-1}, \ldots, t^1, 0; p) : (t^1, \ldots, t^{s-1}, t^s) \in \mathbb{C}^{sn}\}$$

equals $N$. As in [BER99b] we introduce the variables $\eta = (\eta^1, \cdots, \eta^s) \in \mathbb{C}^{sn}$ and $\xi = (\xi^0, \ldots, \xi^{s-1}) \in \mathbb{C}^{sn}$ and put

$$\widetilde{V}(\eta, \xi; x) := v^{2s}(\xi^0, \eta^1, \cdots, \eta^{s-1}, \eta^s, \eta^{s-1} + \xi^1, \eta^{s-2} + \xi^2 \cdots, \eta^1 + \xi^{s-1}; x). \tag{3.5}$$

Then the generic rank of $\frac{\partial \widetilde{V}}{\partial \xi}$ along $\{(\eta, 0; x)\}$ equals $N$ and we can reorder the $sn$ components of $\xi$ and write $\xi = (\xi^1, \xi^2) \in \mathbb{C}^N \times \mathbb{C}^{sn-N}$ such that the determinant of $\frac{\partial \widetilde{V}}{\partial \xi^1}(\eta, 0; x)$ does not vanish identically. Hence the map

$$V(\eta, \xi^1; x) := \widetilde{V}(\eta, \xi^1, 0; x) \tag{3.6}$$

satisfies the assumptions of the following lemma with $\mu := \kappa - k$ (where $\xi^1$ is replaced by $\xi$ for brevity). The lemma is a variant of Proposition 4.1.18 in [BER99b] with smooth parameters (see also Lemma 10.5 in [BRZ00]). The proof closely follows the one given in [BER99b].

**Lemma 3.3.** *Let $V(\eta, \xi; x)$ be a map of class $\mathcal{C}^\mu$ from a neighborhood of $0$ in $\mathbb{C}^h \times \mathbb{C}^N \times \mathbb{R}^m$ into $\mathbb{C}^N$ with $V(\cdot, \cdot; x)$ holomorphic for each fixed $x$. Suppose that $V(\eta, 0; x) \equiv 0$ and*

$$\delta(\eta; x) := \det\left(\frac{\partial V}{\partial \xi}(\eta, 0; x)\right)^2 \not\equiv 0 \tag{3.7}$$

*for each fixed $x$. Then there exists a map $\varphi(\eta, \widetilde{Z}; x)$ from a neighborhood of $0$ in $\mathbb{C}^h \times \mathbb{C}^N \times \mathbb{R}^m$ into $\mathbb{C}^N$ of class $\mathcal{C}^\mu$ and holomorphic for each fixed $x$, such that the identity*

$$V\left(\eta, \varphi\left(\eta, \frac{Z}{\delta(\eta; x)}; x\right); x\right) \equiv Z$$

*holds for all $(\eta, Z; x)$ with $(\eta; x)$ and $\frac{Z}{\delta(\eta; x)}$ sufficiently small.*

*Proof.* Since $V(\eta, \xi; x)$ is holomorphic in $(\eta, \xi)$ and $V(\eta, 0; x) \equiv 0$, we can write it in the form

$$V(\eta, \xi; x) = V_1(\eta, \xi; x)\xi,$$

where $V_1$ is the $N \times N$ matrix function of class $\mathcal{C}^\mu$ given by

$$V_1(\eta, \xi; x) := \int_0^1 \frac{\partial V}{\partial \xi}(\eta, t\xi; x)dt.$$

In particular, we have

$$V_1(\eta, 0; x) \equiv \frac{\partial V}{\partial \xi}(\eta, 0; x). \tag{3.8}$$

The map $V_1$ is also of class $\mathcal{C}^\mu$ because of the Cauchy type formula for the derivatives. By applying the same procedure one more time we can write $V_1(\eta, \xi; x) = V_1(\eta, 0; x) + V_2(\eta, \xi; x)\xi$ with $V_2$ again being of class $\mathcal{C}^\mu$. Then

$$V(\eta, \xi; x) = V_1(\eta, 0; x)\xi + R(\eta, \xi; x)(\xi, \xi),$$

where $R(\eta, \xi; x)(\cdot, \cdot)$ is the bilinear form defined by $V_2$.

By (3.7) and (3.8) we have $\det(V_1(\eta, 0; x)) \not\equiv 0$ for each fixed $x$. We set $\Delta(\eta; x) := \det(V_1(\eta, 0; x))$. We now proceed exactly as in the proof of Proposition 4.1.18 in [BER99b]. We write the equation $V(z, \xi; x) = Z$ as

$$V_1(\eta, 0; x)\xi + R(\eta, \xi; x)(\xi, \xi) = Z$$

and solve it for $\xi$ by using Cramer's rule in terms of $(Z)/\Delta(\eta; x)$ and $R(\eta, \xi; x)(\xi/\Delta(\eta; x), \xi)$. In the expression obtained we divide both sides by $\Delta(\eta; x)$ and apply the implicit function theorem with respect to $\xi/\Delta(\eta; x)$:

$$\frac{\xi}{\Delta(\eta; x)} = \psi\Big(\eta, \frac{Z}{\Delta(\eta; x)^2}; x\Big), \quad \psi \in \mathcal{C}^\mu, \quad \psi(\cdot, \cdot; x) \text{ holomorphic.}$$

Then the map

$$\varphi(\eta, \widetilde{Z}; x) := \Delta(\eta; x)\psi(\eta, \widetilde{Z}; x) \in \mathbb{C}^N$$

satisfies the conclusion of the lemma. $\qquad\blacksquare$

### 3.4. Hörmander numbers and the vanishing order of $\delta(\eta; x)$.

As before we assume that $M$ has the finite type $\nu \leq \kappa$ at $p$. Recall that the *Hörmander numbers* $2 \leq \mu_1 \leq \ldots \leq \mu_d = \nu$ *of* $M$ *at* $p$ are defined as follows. For $1 \leq j \leq d$, the $j$th Hörmander number $\mu_j$ is the minimal number $\mu \leq \nu$ with the property that all vector fields in $T^c M$ (of class $\mathcal{C}^{\kappa-1}$) together with their commutators of length less or equal to $\mu$ span at $p$ a vector subspace of $T_p M$ of dimension at least $n + j$. In particular, we have $\mu_d = \nu$. (This definition differs slightly from the one given in §3.4 in [BER99a], where the same Hörmander numbers are counted as a single number with multiplicity).

As before we choose a family of approximate embeddings $M(x) \subset \mathbb{C}^N$ of order $k$ with $\nu \leq k \leq \kappa$. Then $M(p)$ is real-analytic and also has the type $\nu$. Moreover, since the values at $p$ of all commutators of vector fields in $T^c M(p)$ of length at most $k$ coincide with those of vector fields in $T^c M$, all Hörmander numbers of $M(p)$ at $p$ coincide with those of $M$. We now use the so-called *normal canonical coordinates* (see [BG77, BR87]). For our purposes, it will be convenient to have the normal canonical coordinates in the form of Theorem 3.2.3 in [BER99b] (see also Theorem 4.5.1 in [BER99a]). It follows that there exist holomorphic local coordinates $Z = (z, w) \in \mathbb{C}^n \times \mathbb{C}^d = \mathbb{C}^N$ in a neighborhood of $p$, where $M(p)$ is locally given by the equations

$$w_j = \overline{w}_j + P_j(z, \overline{z}, \overline{w}) + R_j(z, \overline{z}, \overline{w}), \quad 1 \leq j \leq d,$$

where the functions $P_j$ and $R_j$ are holomorphic in their arguments, satisfy the normalization conditions

$$P_j(z, 0, \overline{w}) \equiv R_j(z, 0, \overline{w}) \equiv P_j(0, \overline{z}, \overline{w}) \equiv R_j(0, \overline{z}, \overline{w}) \equiv 0$$

and, in the notation $\varepsilon^\mu w := (\varepsilon^{\mu_1} w_1, \ldots, \varepsilon^{\mu_d} w_d)$, the weight homogeneity conditions

$$P_j(\varepsilon z, \varepsilon \overline{z}, \varepsilon^\mu \overline{w}) = \varepsilon^{\mu_j} P_j(z, \overline{z}, \overline{w}), \quad R_j(\varepsilon z, \varepsilon \overline{z}, \varepsilon^\mu \overline{w}) = O(\varepsilon^{\mu_j + 1})$$

for small $\varepsilon > 0$. The functions $P_j(z, \overline{z}, \overline{w})$ are nonzero polynomials and the submanifold $M_0 := \{w_j = \overline{w}_j + P_j(z, \overline{z}, \overline{w}), 1 \leq j \leq d\}$ is of the same type $\nu$ and has the same Hörmander numbers ($M_0$ is a *weighted homogeneous generic submanifold* in the sense of §4.4 in [BER99a]).

As was observed in [BER99b], in the canonical coordinates, the map $v^s(\cdot, \cdot; p)$ has a power series expansion such that its first $n$ components are linear and, for $1 \leq j \leq d$, its $(n+j)$th component starts from a nonzero homogeneous polynomial of degree $\mu_j$. The same holds for the components of $V$ (defined by (3.5) and (3.6)). Hence the power series expansion of the determinant in (3.7) starts from a homogeneous polynomial of degree $(\mu_1 - 1) + \ldots + (\mu_d - 1) = \mu_1 + \cdots + \mu_d - d$. Since $M_0$ is of finite type at $p$, the determinant cannot vanish identically. Thus we come to the following conclusion:

**Lemma 3.4.** *Let $V(\eta, \xi; x)$ be as before and $\delta(\eta; x)$ be given by Lemma 3.3. Then there exists a vector $\eta_0 \in \mathbb{C}^{sn}$ such that the vanishing order $m$ of the function $\widetilde{\delta}(\lambda) := \delta(\lambda \eta_0; p)$ at $\lambda = 0 \in \mathbb{C}$ equals $2(\mu_1 + \cdots + \mu_d - d)$. In particular, $m \leq 2d(\nu - 1)$.*

## 4. Jet parametrization of local holomorphic maps

Consider two embedded CR-submanifolds $M \subset \mathbb{C}^N$ and $M' \subset \mathbb{C}^{N'}$. If $p \in M$ and $p' \in M'$ are fixed and $f: (M, p) \to (M', p')$ is a germ of a $\mathcal{C}^r$-smooth CR-map, it can be shown (see e.g. the arguments in the proof of Proposition 1.7.14 in [BER99a]) that $f$ can be approximated up to order $r$ by the restriction to $M$ of a holomorphic $r$th order polynomial map $F: \mathbb{C}^N \to \mathbb{C}^{N'}$ satisfying

$$\rho'(F(x), \overline{F(x)}) = o(|x - p|^r) \text{ as } x \to p \text{ in } M. \tag{4.1}$$

If (4.1) holds, we say that $F$ *sends $M$ into $M'$ up to order $r$ at $p$*. This notion generalizes that of $(r + 1)$-equivalence given in [BRZ00]. If $M$ is generic at $p$, i.e. if $T_p M + J T_p M = T_p \mathbb{C}^N$ for any $p \in M$, then the jet $j_p^r F \in J^r(\mathbb{C}^N, \mathbb{C}^{N'})$ is uniquely determined by $j_p^r f \in J^r(M, M')$. The condition (2.1) of $l$-nondegeneracy can be also written with $f$ replaced by $F$ and is again independent of the choice of $L_1, \ldots, L_n$ and of $Z'$ provided $F$ is holomorphic. If, in addition, $F$ sends $M$ into $M'$ up to order $r$ at $x$, then this condition is also independent of the choice of $\rho'$. In this case we call both $F$ and its jet $j_p^r F$ $l$-nondegenerate with respect to $M$ and $M'$. We denote by $D_{M, M'}^{r, l}(\mathbb{C}^N, \mathbb{C}^{N'}) \subset J^r(\mathbb{C}^N, \mathbb{C}^{N'})$ the open subset of jets whose representatives either do not send $M$ into $M'$ up to order $r$ or are $l$-nondegenerate with respect to $M$ and $M'$.

The following is the main technical result of this section:

**Theorem 4.1.** *Let $M \subset \mathbb{C}^N$ and $M' \subset \mathbb{C}^{N'}$ be $\mathcal{C}^\kappa$-smooth generic submanifolds ($2 \leq \kappa \leq \infty$). Suppose that $M$ is of codimension $d \geq 1$ and type $2 \leq \nu \leq \kappa$. Fix a number $l \geq 0$, set $r := 2(d + 1)l$ and fix a jet $(x_0, \Lambda_0) \in D_{M, M'}^{r, l}(\mathbb{C}^N, \mathbb{C}^{N'})$. Then there exist a neighborhood $\Omega =$* $\Omega(x_0, (x_0, \Lambda_0)) \subset \mathbb{C}^N \times J^r(\mathbb{C}^N, \mathbb{C}^{N'})$ and, for every $k$ with $\max(r, \nu) \leq k \leq \kappa - 1$, a $\mathcal{C}^{\kappa - k - 1}$-smooth map $\Psi^k(Z, x, \Lambda) \colon \Omega \to \mathbb{C}^{N'}$, holomorphic in $(Z, \Lambda) \in \mathbb{C}^N \times J^r_{x,x'}(\mathbb{C}^N, \mathbb{C}^{N'})$ for each fixed $(x, x') \in M \times M'$, such that the following holds. For every $x \in M$, every holomorphic map $F$ from a neighborhood of $x$ in $\mathbb{C}^N$ into $\mathbb{C}^{N'}$ with $(x, x, j^r_x F) \in \Omega$ that sends $M$ into $M'$ up to order $k$ at $x$, satisfies the identity

$$F(Z) = \Psi^k(Z, x, j^r_x F) + o\big(|Z - x|^{\frac{k-r}{m+1}}\big), \quad Z \to x, \tag{4.2}$$

where $m$ is the vanishing order of the function $\delta(\cdot, p)$ given by (3.7). If $M$ and $M'$ are real-analytic, the maps $\Psi^k$ can be also chosen to be real-analytic. Moreover, if $M$ and $M'$ depend in a $\mathcal{C}^\kappa$ fashion (resp. real-analytically) on some real parameters, then the maps $\Psi^k$ can be also chosen to depend in a $\mathcal{C}^{\kappa - k - 1}$ fashion (resp. real-analytically) on the parameters.

The rest of this section is devoted to the proof of Theorem 4.1 that is splitted into three steps that may be of independent interest.

4.1. **Basic reflection identity.** Our first goal will be to establish a relation between jets of a holomorphic map $F$ as in Theorem 4.1 at two different points $Z$ and $\zeta$ satisfying the reflection relation $(Z, \zeta) \in \mathcal{M}(x)$, where $\mathcal{M}(x)$ is the approximate complexification of $M$ as constructed in §3.2. Let $\rho$ be a defining function for $M$ near a point $p \in M$. For $l \leq k \leq \kappa$ and $x \in M$ near $p$, consider approximations $M(x)$ with real-analytic defining functions $\rho(Z, \overline{Z}; x)$ as defined in §3.2. In the sequel we shall find it convenient to work with $(1, 0)$ vector fields on $M$ that are precisely the conjugates of the $(0, 1)$ vector fields (see §1). We write $L_j = L_j(Z, \overline{Z})$, $1 \leq j \leq n$, for a basis of $\mathcal{C}^{\kappa - 1}$-smooth $(1, 0)$ vector fields on $M$ near $p$ and, for each $x \in M$ near $p$, approximate them by restrictions of a basis $L_j(x) = L_j(Z, \overline{Z}; x)$ of real-analytic $(1, 0)$ vector fields in $\mathbb{C}^N$ such that

$$L(x)^\alpha \chi(Z)|_{Z=x} = L^\alpha \chi(Z)|_{Z=x}, \quad x \in M,$$

for any smooth function $\chi$ in a neighborhood of $x$ and for all multiindices $\alpha$ with $|\alpha| \leq l$. This means that the restrictions to $M$ of the coefficients of $L_j(x)$ have the same $(l - 1)$-jets at $x$ as the coefficients of $L_j$. The approximations $L_j(x) = L_j(Z, \overline{Z}; x)$ can be chosen to be of class $\mathcal{C}^{\kappa - l}$ in all their variables including $x$.

For the proof of Theorem 4.1, we fix an $l$-nondegenerate holomorphic $r$-jet

$$(x_0, \Lambda_0) \in D^{r,l}_{M,M'}(\mathbb{C}^N, \mathbb{C}^{N'})$$

sending $M$ into $M'$ and $x_0$ to a fixed point $x'_0 \in M'$ and consider a point $x \in M$ and a holomorphic map $F$ from a neighborhood of $x$ in $\mathbb{C}^N$ into $\mathbb{C}^{N'}$ sending $M$ into $M'$ up to order $k$ at $x$. We shall assume the $r$-jet $j^r_x F$, $x \in M$, to be sufficiently close to $(x_0, \Lambda_0)$. Then, for $f := F|_M$, the condition (2.1) still holds with $p$ replaced by $x$. Moreover, the vector fields $L_1, \ldots, L_n$ in (2.1) can be replaced by their approximations $L_1(x), \ldots, L_n(x)$ without changing the span in (2.1). We conclude that

$$\mathrm{span}\{L(x)^\alpha \rho'^j_{Z'}(F(Z), \overline{F(Z)})(x, \overline{x}) : 1 \leq j \leq d', 0 \leq |\alpha| \leq l\} = \mathbb{C}^{N'} \tag{4.3}$$

for $F$ as above.

We use the notation from §3 for $M$ and analogous notation for $M'$. Let $M'(x') \subset \mathbb{C}^{N'}$ be a family of approximations for $M'$ with defining functions $\rho'(Z', \overline{Z'}; x')$. Then since $F$ sends $M$ into $M'$ up to order $k$ at $x$, we conclude

$$\rho'(F(Z), \overline{F(Z)}; F(x)) = o(|Z - x|^k) \text{ as } Z \to x \text{ in } M(x). \tag{4.4}$$

By applying the vector fields $L_j(x)$ to (4.4) we obtain

$$L(x)^\alpha \rho'(F(Z), \overline{F(Z)}; F(x)) = o(|Z - x|^{k-|\alpha|}) \text{ as } Z \to x \text{ in } M(x) \tag{4.5}$$

for $|\alpha| \leq l$. Since each $L_j(Z, \overline{Z}; x)$ is real-analytic in $(Z, \overline{Z})$, we can consider the complexifications $L_j(Z, \zeta; x)$. By the standard complexification argument, (4.5) implies

$$L(Z, \zeta; x)^\alpha \rho'(F(Z), \overline{F}(\zeta); F(x)) = o(|(Z - x, \zeta - \overline{x})|^{k-|\alpha|}) \tag{4.6}$$

as $(Z, \zeta) \to (x, \overline{x})$ in $\mathcal{M}(x)$. Recall that $\rho'(Z', \zeta'; x')$ is of class $\mathcal{C}^{\kappa-k}$ and each $L(Z, \zeta; x)$ is of class $\mathcal{C}^{\kappa-l}$ and they are holomorphic in $(Z', \zeta')$ and in $(Z, \zeta)$ respectively. Under the assumptions of Theorem 4.1 we clearly have $l \leq k$ and hence $L(Z, \zeta; x)$ is in $\mathcal{C}^{\kappa-k}$. Since all differentiations in (4.6) are taken in the "holomorphic directions", we conclude that the left-hand side of (4.6) is of class $\mathcal{C}^{\kappa-k}$. By using the chain rule and the fact that $L(Z, \zeta; x)^\alpha \overline{F}(\zeta) \equiv 0$, we can write the left-hand side of (4.6) for $|\alpha| \leq l$ as a $\mathcal{C}^{\kappa-k}$-smooth function in $\left(Z, j_Z^l F, \zeta, \overline{F}(\zeta); x, F(x)\right)$, holomorphic in $\left(Z, j_Z^l F, \zeta, \overline{F}(\zeta)\right)$ for each fixed $(x, F(x))$.

By (4.3), the implicit function theorem can be used to solve (4.6) for $\overline{F}(\zeta)$ in the form

$$\overline{F}(\zeta) = \Psi^0(\zeta, Z, j_Z^l F; x, F(x)) + o(|(Z - x, \zeta - \overline{x})|^{k-l}) \tag{4.7}$$

as $(Z, \zeta) \to (x, \overline{x})$ in $\mathcal{M}(x)$, where $\Psi^0(\zeta, Z, \Lambda; x, x')$ is a $\mathbb{C}^{N'}$-valued function of class $\mathcal{C}^{\kappa-k}$ in a neighborhood of $(\overline{x_0}, x_0, \Lambda_0; x_0, x'_0)$, holomorphic in $(\zeta, Z, \Lambda)$ for each fixed $(x, x') \in M \times M'$.

The identity (4.7) holds for $(Z, \zeta) \in \mathcal{M}(x)$. Due to the construction of $\mathcal{M}(x)$, for each parameter $(Z_0, \zeta_0; x) \in \mathbb{C}^N \times \mathbb{C}^N \times M$ close to $(x_0, \overline{x_0}; x_0)$, there exists a holomorphic function $Z = Z(\zeta) \in \mathbb{C}^N$ in a neighborhood of $\zeta_0$ with $Z(\zeta_0) = Z_0$, depending in a $\mathcal{C}^{\kappa-k}$ fashion on all variables and parameters and satisfying $(Z(\zeta), \zeta) \in \mathcal{M}(x)$. After substituting $Z(\zeta)$ for $Z$ in (4.7) and differentiating in $\zeta$ we obtain the *basic reflection identity*

$$j_\zeta^\tau \overline{F} = \Psi^\tau(\zeta, Z, j_Z^{\tau+l} F; x, F(x)) + o(|(Z - x, \zeta - \overline{x})|^{k-l-\tau}) \tag{4.8}$$

as $(Z, \zeta) \to (x, \overline{x})$ in $\mathcal{M}(x)$, $\tau \leq k - l$, where $\Psi^\tau(\zeta, Z, \Lambda; x, x')$ is still of class $\mathcal{C}^{\kappa-k}$, because we differentiated only in the "holomorphic direction".

Observe that the $r$-jet $\Lambda_0 \in J_{x_0}^r(\mathbb{C}^N, \mathbb{C}^{N'})$ is assumed to send $M$ into $M'$. Then, by following the construction, we see that (4.8) also holds with $k$ replaced by $r$ and any representative $F$ of $(x_0, \Lambda_0)$. Hence we obtain

$$\overline{\Lambda_0^\tau} = \Psi^\tau(\overline{x_0}, x_0, \Lambda_0^{\tau+l}; x_0, x'_0) \tag{4.9}$$

for $\tau + l \leq r$, where $\Lambda_0^\tau \in J^\tau(\mathbb{C}^N, \mathbb{C}^{N'})$ denotes image of $\Lambda_0$ under the canonical projection $J^r(\mathbb{C}^N, \mathbb{C}^{N'}) \to J^\tau(\mathbb{C}^N, \mathbb{C}^{N'})$.

We summarize the result of this paragraph indicating the dependence of $\Psi^\tau$ on $k$ by adding a superscript:

**Proposition 4.2.** *Let $M \subset \mathbb{C}^N$ and $M' \subset \mathbb{C}^{N'}$ be $\mathcal{C}^\kappa$-smooth generic submanifolds ($2 \leq \kappa \leq \infty$). Fix any numbers $l, \tau, k \geq 0$ with $\tau + l \leq k \leq \kappa - 1$ and a jet $(x_0, \Lambda_0) \in D^{\tau + l, l}_{M,M'}(\mathbb{C}^N, \mathbb{C}^{N'})$ with target point $x'_0 \in M'$ that sends $M$ into $M'$. Let $M(x) \subset \mathbb{C}^N$ be any family of approximate real-analytic embeddings of $M$ of order $k$ as in §3.2 for $x \in M$ near $x_0$. Then there exist a neighborhood*

$$\widetilde{\Omega} = \widetilde{\Omega}(\overline{x}_0, (x_0, \Lambda_0); x_0, x'_0) \subset \mathbb{C}^N \times J^{\tau + l}(\mathbb{C}^N, \mathbb{C}^{N'}) \times M \times M'$$

*and a $\mathcal{C}^{\kappa - k}$-smooth map*

$$\Psi^{k,\tau}(\zeta, Z, \Lambda; x, x') \colon \widetilde{\Omega} \to J^\tau(\mathbb{C}^N, \mathbb{C}^{N'}),$$

*holomorphic in $(\zeta, Z, \Lambda) \in \mathbb{C}^N \times J^r_{x,x'}(\mathbb{C}^N, \mathbb{C}^{N'})$ for each fixed $(x, x') \in M \times M'$, such that the following holds.*

(i) *For the image $\Lambda_0^\tau \in J^\tau(\mathbb{C}^N, \mathbb{C}^{N'})$ of $\Lambda_0$ under the projection $J^{\tau + l}(\mathbb{C}^N, \mathbb{C}^{N'}) \to J^\tau(\mathbb{C}^N, \mathbb{C}^{N'})$, one has $\overline{\Lambda_0^\tau} = \Psi^{k,\tau}(\overline{x_0}, x_0, \Lambda_0; x_0, x'_0)$.*

(ii) *For every $x \in M$, every holomorphic map $F$ from a neighborhood of $x$ in $\mathbb{C}^N$ into $\mathbb{C}^{N'}$ with $(\overline{x}, x, j^r_x F; x, F(x)) \in \widetilde{\Omega}$ that sends $M$ into $M'$ up to order $k$ at $x$, satisfies the identity*

$$j^\tau_\zeta \overline{F} = \Psi^{k,\tau}(\zeta, Z, j^{\tau + l}_Z F; x, F(x)) + o(|(Z - x, \zeta - \overline{x})|^{k - l - \tau}) \tag{4.10}$$

*as $(Z, \zeta) \to (x, \overline{x})$ in $\mathcal{M}(x)$.*

*If $M$ and $M'$ are real-analytic, the family $M(x)$ and the maps $\Psi^{k,\tau}$ can be also chosen to be real-analytic. Moreover, if $M$ and $M'$ depend in a $\mathcal{C}^\kappa$ fashion (resp. real-analytically) on some real parameters, then the maps $\Psi^{k,\tau}$ can be also chosen to depend in a $\mathcal{C}^{\kappa - k}$ fashion (resp. real-analytically) on the parameters.*

**4.2. Singular jet parametrization.** Our next step is to iterate the identity (4.10) as in [Z97, BER99b, BRZ00]. The property (i) in Proposition 4.2 allows us to iterate it for $(x, j^r_x F)$ sufficiently close to $(x_0, \Lambda_0)$. After iterating (4.10) $2(d+1)$ times we obtain

$$F(\zeta^0) = \Phi^k(\zeta^0, \dots, \zeta^{2(d+1)}, j^r_{\zeta^{2(d+1)}} F; x, F(x)) +$$

$$o(|(\zeta^0 - x, \zeta^1 - \overline{x}, \dots, \zeta^{2(d+1)} - x)|^{k-r}) \tag{4.11}$$

as $(\zeta^0, \dots, \zeta^{2(d+1)}) \to (x, \overline{x}, \dots, x)$ in $\mathcal{M}^{2(d+1)}$, where $\Phi^k(\zeta^0, \dots, \zeta^{2(d+1)}, \Lambda; x, x')$ is of class $\mathcal{C}^{\kappa - k}$ and holomorphic in $(\zeta^0, \dots, \zeta^{2(d+1)}, \Lambda)$.

Now let $\delta$ and $\varphi$ be as in Lemma 3.3. Put $\zeta^{2(d+1)} = x$ and

$$\zeta^0 = Z = V\left(\eta, \varphi\left(\eta, \frac{Z - x}{\delta(\eta; x)}; x\right); x\right) + x$$

in the identity (4.11). We obtain

$$F(Z) \equiv \widetilde{\Psi}^k\left(\eta, \frac{Z - x}{\delta(\eta; x)}, j^r_x F; x, F(x)\right) + R^k_F\left(\eta, \frac{Z - x}{\delta(\eta; x)}; x\right), \tag{4.12}$$

where $\widetilde{\Psi}^k$ and $R^k_F$ are of class $\mathcal{C}^{\kappa - k}$, their restrictions $\widetilde{\Psi}^k(\cdot, \cdot, \cdot; x, x')$ and $R^k_F(\cdot; x)$ are holomorphic and

$$R^k_F(\eta, \widetilde{Z}; x) = o(|(\eta, \widetilde{Z})|^{k-r}) \text{ as } (\eta, \widetilde{Z}) \to 0 \tag{4.13}$$

for each fixed $x$. We note that here, exactly as in [BRZ00], only $R_F^s$ (but not $\widetilde{\Psi}^s$) depends on $F$ (this dependence is indicated by $F$ in the subscript). A new ingredient comparing with [BRZ00] is that $\widetilde{\Psi}^k$ depends on the approximations $M(x)$ and $M'(x')$, in particular, also on $k$.

Let $\eta_0 \in \mathbb{C}^{sn}$ be given by Lemma 3.4 and define $\widehat{\delta}(\lambda; x) := \delta(\lambda \eta_0; x)$ (for $\lambda \in \mathbb{C}$). Then, setting $\eta = \lambda \eta_0$ in (4.12), we have

$$F(Z) \equiv \widehat{\Psi}^k\left(\lambda, \frac{Z-x}{\widehat{\delta}(\lambda; x)}, j_x^r F; x, F(x)\right) + \widehat{R}_F^k\left(\lambda, \frac{Z-x}{\widehat{\delta}(\lambda; x)}; x\right), \tag{4.14}$$

where $\widehat{\Psi}^k$ and $\widehat{R}_F^k$ are defined in the obvious way and have the properties similar to those of $\widetilde{\Psi}^k$ and $R_F^k$ respectively. After a holomorphic change of the parameter $\lambda \in \mathbb{C}$, we may assume that $\widehat{\delta}(\lambda; x_0) = \lambda^m$, where $m \leq 2d(\nu-1)$ is the vanishing order considered in Lemma 3.4. We choose a smooth parametrization $s \mapsto x(s)$ of $M$ near $x_0$ sending 0 to $x_0$. We view $s$ as a row in $\mathbb{R}^{\dim M}$. Then we have

$$\widehat{\delta}(\lambda; x(s)) \equiv \lambda^m + a(\lambda; s)s, \tag{4.15}$$

where $a(\lambda; s)$ is a $\mathcal{C}^{\kappa-k-1}$-smooth column which is holomorphic for each fixed $s$. Then (4.14) implies

$$F(Z) \equiv \widehat{\Psi}^k\left(\lambda, \frac{Z-x(s)}{\lambda^m}\left(1+\frac{a(\lambda; s)s}{\lambda^m}\right)^{-1}, j_{x(s)}^r F; x(s), F(x(s))\right)$$
$$+ \widehat{R}_F^k\left(\lambda, \frac{Z-x(s)}{\widehat{\delta}(\lambda; x(s))}; x(s)\right). \tag{4.16}$$

We can now rewrite (4.16) in the form

$$F(Z) \equiv \check{\Psi}^k\left(\lambda, \frac{Z-x(s)}{\lambda^m}, \frac{s}{\lambda^m}, j_{x(s)}^r F; x(s), F(x(s))\right) + \widehat{R}_F^k\left(\lambda, \frac{Z-x(s)}{\widehat{\delta}(\lambda; x(s))}; x(s)\right), \tag{4.17}$$

where $\check{\Psi}^k(\lambda, \widetilde{Z}, \widetilde{x}, \Lambda; x, x')$ is defined and $\mathcal{C}^{\kappa-k-1}$-smooth for $(\lambda, \widetilde{Z}, \widetilde{x}, \Lambda; x, x')$ in a neighborhood of $(0, 0, 0, \Lambda_0; x_0, x_0')$ and holomorphic in $(\lambda, \widetilde{Z}, \widetilde{x}, \Lambda)$. This is the *singular jet parametrization* that we shall use for the proof of Theorem 4.1.

### 4.3. Resolution of singularities with smooth parameters.
We will apply the following statement to the map $\check{\Psi}^k$ that can be seen as a variant of Lemma 10.6 in [BRZ00] with smooth parameters.

**Lemma 4.3.** *Let $P(\lambda, \widetilde{t}; x)$ be a function of class $\mathcal{C}^\tau$ in a neighborhood of 0 in $\mathbb{C} \times \mathbb{C}^n \times \mathbb{R}^h$ such that $P(\cdot, \cdot; x)$ is holomorphic for each fixed $x$. Let*

$$P\left(\lambda, \frac{t}{\lambda^m}; x\right) \equiv \sum_\nu c_\nu(t; x)\lambda^\nu \tag{4.18}$$

*be the Laurent series expansion in $\lambda$. Then the coefficient $c_0$ is of class $\mathcal{C}^\tau$ and is holomorphic in $t$ for each fixed $x$.*

*Proof.* Since $P(\cdot, \cdot; x)$ is holomorphic, we can write

$$P\left(\lambda, \widetilde{t}; x\right) \equiv \sum_{\nu, \alpha} P_{\nu, \alpha}(x) \lambda^{\nu} \widetilde{t}^{\alpha}$$

in a neighborhood of the origin. Then $P_{\nu, \alpha} \in \mathcal{C}^{\tau}$ and

$$c_0(t; x) \equiv \sum_{\nu = |\alpha| m} P_{\nu, \alpha}(x) t^{\alpha}$$

is clearly holomorphic in $t$ near the origin. Since the partial differentiations commute, all partial derivatives of $P(\lambda, \widetilde{t}; x)$ with respect to $x$ are holomorphic in the variables $(\lambda, \widetilde{t})$:

$$\overline{\partial}_{(\lambda, \widetilde{t})} \partial_x^{\beta} P \equiv \partial_x^{\beta} \overline{\partial}_{(\lambda, \widetilde{t})} P \equiv 0$$

for every $|\beta| \leq \tau$. Clearly the coefficient $c_0(t; x)$ does not change if $\lambda$ is replaced by $C\lambda$ in (4.18) with $C \neq 0$ being any constant. Hence we may assume that $P$ and all partial derivatives in $x$ up to order $\tau$ are holomorphic for $|\lambda| < 2$, $|\widetilde{t}| < \varepsilon$. This means that

$$\partial_x^{\beta} P(\lambda, t; x) \equiv \sum_{\nu, \alpha} \partial_x^{\beta} P_{\nu, \alpha}(x) \lambda^{\nu} t^{\alpha}$$

converges absolutely on $|\lambda| \leq 1, |\widetilde{t}| < \varepsilon$ and so does

$$\partial_x^{\beta} c_0(t; x) \equiv \sum_{\nu = |\alpha| m} \partial_x^{\beta} P_{\nu, \alpha}(x) t^{\alpha}.$$

This shows that $c_0(t; x)$ is of class $\mathcal{C}^{\tau}$ as required.                     □

The following lemma will be used for $\widehat{R}_F^k$ and can be also seen as a variant of Lemma 10.6 in [BRZ00] with smooth parameters.

**Lemma 4.4.** *Let $R(\lambda, \widetilde{t})$ and $\delta(\lambda)$ be holomorphic functions in $\mathbb{C} \times \mathbb{C}^n$ and in $\mathbb{C}$ respectively near the origins such that*

$$R(\lambda, \widetilde{t}) = o(|(\lambda, \widetilde{t})|^h) \text{ as } (\lambda, \widetilde{t}) \to 0$$

*and suppose that $\delta(\lambda)$ has the vanishing order $g$ at $\lambda = 0$ for some nonnegative integers $h$ and $g$. Consider the Laurent series expansion*

$$R\left(\lambda, \frac{t}{\delta(\lambda)}\right) = \sum_{\nu} c_{\nu}(t) \lambda^{\nu}.$$

*Then $c_0(t) = o(|t|^{\frac{h}{g+1}})$ as $t \to 0$.*

*Proof.* We write $R$ and $\frac{t}{\delta(\lambda)}$ as power series:

$$R(\lambda, \widetilde{t}) \equiv \sum_{\nu, \alpha} R_{\nu, \alpha} \lambda^{\nu} \widetilde{t}^{\alpha}, \quad \frac{t}{\delta(\lambda)} \equiv t\lambda^{-g}(a_0 + a_1 \lambda + a_2 \lambda^2 + \dots),$$

where $a_0 \neq 0$. The substitution $\widetilde{t} = \frac{t}{\delta(\lambda)}$ yields

$$R\Big(\lambda, \frac{t}{\delta(\lambda)}\Big) \equiv \sum_{\nu,\alpha} R_{\nu,\alpha} t^\alpha \lambda^{\nu - g|\alpha|} (a_0 + a_1\lambda + a_2\lambda^2 + \dots)^\alpha.$$

Then it is easy to see that $c_0(t)$ is a convergent power series whose entries are linear combinations of the terms $R_{\nu,\alpha} t^\alpha$ with $\nu - g|\alpha| \leq 0$. By the assumption, $R(\lambda, \widetilde{t}) = o(h)$ and hence the inequality $\nu + |\alpha| > h$ holds for all nonvanishing coefficients $R_{\nu,\alpha}$. Putting these inequalities together we obtain $|\alpha| > \frac{h}{g+1}$ showing $c_0(t) = o(\frac{h}{g+1})$ as required. $\qquad\square$

### 4.4. The end of proof of Theorem 4.1.
The first statement of Theorem 4.1 is a consequence of Proposition 4.2 and the following statement:

**Proposition 4.5.** *Let $M, M', \kappa, \nu, d, l, r$ be as in Theorem 4.1, fix a jet $(x_0, \Lambda_0) \in J^r_{M,M'}(\mathbb{C}^N, \mathbb{C}^{N'})$ that sends $M$ into $M'$ and let $M(x)$ be as in Proposition 4.2. For each $0 \leq \tau \leq 2dl$, denote by $\Lambda_0^{\tau+l} \in J^\tau(\mathbb{C}^N, \mathbb{C}^{N'})$ the image of $\Lambda_0$ under the projection $J^r(\mathbb{C}^N, \mathbb{C}^{N'}) \to J^{\tau+l}(\mathbb{C}^N, \mathbb{C}^{N'})$ and let further, as in Proposition 4.2,*

$$\Psi^{k,\tau} \colon \widetilde{\Omega}^\tau\big(\overline{x}_0, (x_0, \Lambda_0^{\tau+l}); x_0, x_0'\big) \to J^\tau(\mathbb{C}^N, \mathbb{C}^{N'})$$

*be a $\mathcal{C}^{\kappa-k}$-smooth map, holomorphic in $(\zeta, Z, \Lambda) \in \mathbb{C}^N \times J^r_{x,x'}(\mathbb{C}^N, \mathbb{C}^{N'})$ for each fixed $(x, x') \in M \times M'$, having property* (i) *of Proposition 4.2. Then there exists a $\mathcal{C}^{\kappa-k-1}$-smooth map $\Psi^k \colon \Omega(x_0, (x_0, \Lambda_0)) \to \mathbb{C}^{N'}$ as in Theorem 4.1 such that the following holds. For every $x \in M$, every holomorphic map $F$ from a neighborhood of $x$ in $\mathbb{C}^N$ into $\mathbb{C}^{N'}$ with $(x, x, j^r_x F) \in \Omega$ satisfying* (4.10) *for all $0 \leq \tau \leq 2dl$, satisfies also the identity* (4.2). *If $M$, $M'$ and $\Psi^{\tau,k}$ are real-analytic, the maps $\Psi^k$ can be also chosen to be real-analytic.*

For the proof, consider the the Laurent series expansion in $\lambda$ for the function $\check{\Psi}^k$ in (4.17):

$$\check{\Psi}^k\Big(\lambda, \frac{\widetilde{Z}}{\lambda^m}, \frac{t}{\lambda^m}, \Lambda; x, x'\Big) \equiv \sum_\nu c^k_\nu(\widetilde{Z}, t, \Lambda; x, x')\lambda^\nu. \tag{4.19}$$

By Lemma 4.3, $c_0(Z, \widetilde{Z}, \Lambda; x, x')$ is of class $\mathcal{C}^{\kappa-k-1}$ and holomorphic in $(\widetilde{Z}, \widetilde{x}, \Lambda)$ for each fixed $(x, x')$. We obtain the $\mathcal{C}^{\kappa-k-1}$-smooth map

$$\Psi^k(\widetilde{Z}, t, \Lambda; x, x') := c^k_0(\widetilde{Z}, t, \Lambda; x, x')$$

defined in a neighborhood of $(x_0, 0, \Lambda_0; x_0, x_0')$ and holomorphic in $(\widetilde{Z}, t, \Lambda)$. The proof of Proposition 4.5 is completed by the following lemma:

**Lemma 4.6.** *Under the assumptions of Theorem 4.1 we have*

$$F(Z) = \Psi^k\big(Z, s, j^r_{x(s)}F; x(s), F(x(s))\big) + o(|Z - x(s)|^{\frac{\kappa-r}{m+1}}) \text{ as } Z \to x(s) \text{ in } \mathbb{C}^N. \tag{4.20}$$

*Proof.* Consider the Laurent series expansions of both sides of (4.17) in $\lambda$. Since $R_F^k = o(k - r)$ by (4.13), it follows from Lemma 4.4 that the constant term in $\lambda$ in the Laurent series expansion of

$$\widehat{R}_F^k\Big(\lambda, \frac{Z - x}{\widehat{\delta}(\lambda; x)}; x\Big) \tag{4.21}$$

is of order $o(\frac{k-r}{m_x+1})$, where $m_x$ is the vanishing order of $\widehat{\delta}(\lambda; x)$ at $\lambda = 0$ for fixed $x$. Clearly we have $m_x \leq m$ for $x \in M$ close to $x_0$, where $m = m_{x_0}$ is the vanishing order of $\widehat{\delta}(\lambda; x_0)$. Therefore the expression in (4.21) is always of order $o(\frac{k-r}{m+1})$. By equating the constant terms (in $\lambda$) in the Laurent series expansion of (4.17) we obtain the required identity. □

It remains to observe that the second statement of Theorem 4.1 about the dependence on a real parameter $t$ can be obtained by adding $t$ to the smooth parameters $x, x'$ in the above proof. The proof of Theorem 4.1 is complete.

## 5. Applications of Theorem 4.1: Proofs of Statements 2.1–2.4

*Proof of Theorem 2.1.* We use Proposition 3.1 for both $M$ and $M'$ to obtain approximate embeddings in $\mathbb{C}^N$ and $\mathbb{C}^{N'}$ parametrized by points $x \in M$ and $x' \in M'$ respectively and then use the statement of Theorem 4.1 for these embeddings. By abusing the notation we write $M$ and $M'$ also for the corresponding embedded generic CR-submanifolds. Choose $(x_0, \lambda_0) \in D^{r,l}(M, M')$ and let $f$ and $x$ satisfy the assumptions of Theorem 2.1. As observed in §4, there exists a $k$th order holomorphic polynomial map $F = F_x \colon \mathbb{C}^N \to \mathbb{C}^{N'}$ whose restriction to $M$ approximates $f$ up to order $k$ at $x$ and such that $F$ sends $M$ into $M'$ up to order $k$ at $x$. Since $M$ is generic in $\mathbb{C}^N$, the jet $j_x^r F$ is uniquely determined by $j_x^r f$, more precisely, $j_x^r F = \Xi_1(x, j_x^r f)$, where $\Xi_1$ is a $\mathcal{C}^{\kappa-r}$-smooth function defined in a neighborhood of $(x_0, \lambda_0)$ in $J^r(M, M')$. By the construction $D^{r,l}(M, M')$, the jet $\Lambda_0 := \Xi_1(x_0, \lambda_0)$ satisfies the assumptions of Theorem 4.1. We may also assume that $F$ satisfies the assumptions of Theorem 4.1.

By Theorem 4.1, we have the identity (4.2) that we differentiate $r + 1$ times in $Z$ and evaluate at $Z = x$. Due to our choice of $k$, the error term in (4.2) does not contribute. Hence we obtain $j_x^{r+1} F = \Psi(x, j_x^r F)$, where $\Psi$ can be chosen to be a local section of the bundle $J^{r+1}(\mathbb{C}^N, \mathbb{C}^{N'}) \to J^r(\mathbb{C}^N, \mathbb{C}^{N'})$ defined in a neighborhood of $(x, j_x^r F)$. Moreover, since $\Psi$ is obtained by differentiating $\Psi^k$ only in $Z$, where $\Psi^k$ is holomorphic, it follows from the Cauchy formula that $\Psi$ has the same smoothness as $\Psi^k$.

Moreover, it follows from the chain rule that $j_x^{r+1} f = \Xi_2(x, j_x^{r+1} F)$, where $\Xi_2$ is a $\mathcal{C}^{\kappa-r-1}$-smooth function in a neighborhood of $(x_0, \Psi(x_0, \Lambda_0))$. By composing $\Psi$ with $\Xi_1$ and $\Xi_2$ we obtain $j_x^{r+1} f = \widetilde{\Phi}(x, j_x^r f)$, where $\widetilde{\Phi}$ is a $\mathcal{C}^{\kappa-r-1}$-smooth $J^{r+1}(M, M')$-valued function defined in a neighborhood of $(x_0, \lambda_0)$. Clearly $\widetilde{\Phi}$ can be assumed to be a local section of the fiber bundle $\pi \colon J^{r+1}(M, M') \to J^r(M, M')$.

In order to obtain a global section $\Phi$ of $\pi$ over $D^{r,l}(M, M')$ as claimed we use the fact that $\pi$ carries a canonical affine bundle structure and, therefore, linear combinations of its sections are well-defined provided the sum of the coefficients is 1. Hence we can use a partition of unity to construct a section $\Phi$ as required. □

*Proof of Corollary 2.2.* Let $\Phi$ be given by Theorem 2.1 and choose local coordinates $x = (x^1, \ldots, x^{\dim M}) \in M$ vanishing at the given point $x_0 \in M$. Then the differential system given by (2.4) implies differential equations for $\Lambda(x) := j_x^r f$ of the form

$$\partial_{x_j} \Lambda(x) = \Phi_j(x, \Lambda(x)). \tag{5.1}$$

Then we construct the required map $\Theta$ by integrating (5.1) along paths. More precisely, every $p$ and $x$ sufficiently close to $x_0$ can be connected in a unique way by a sequence of paths $\gamma_1, \ldots, \gamma_{\dim M}$ such that each $\gamma_j$ is an integral curve of $\partial/\partial x_j$. Then we define $\Theta(x, p, \lambda)$ to be the result of integration (5.1) along this uniquely defined chain of paths. Here, by integration, we mean solving the initial value problem for the corresponding ODE. It is easy to see that $\Theta$ satisfies the required conclusion. $\qquad\square$

*Proof of Corollary 2.4.* Given local coordinates $x = (x_1, \ldots, x_{\dim M}) \in M$, $x' = (x'_1, \ldots, x'_{\dim M'}) \in M'$, we use the associated coordinates $\Lambda = (\Lambda_i^\alpha)$, $1 \le i \le \dim M'$, $\alpha \in \mathbb{Z}_+^{\dim M}$, $0 \le |\alpha| \le r$, for the jet bundle $J_x^r(M, M')$ and write $(x, \Lambda) \in J^r(M, M')$. Let

$$\Phi \colon D^{r,l}(M, M') \subset J^r(M, M') \to J^{r+1}(M, M')$$

be the section given by Theorem 2.1. Then $\Phi$ defines in a natural way a connection on the fibration $\pi \colon D^{r,l}(M, M') \to M$, i.e. a distribution $H$ of horisontal vector subspaces $H(x, \Lambda) \subset T_{(x,\Lambda)} D^{r,l}(M, M')$ such that $d\pi$ restricts to isomorphisms between $H(x, \Lambda)$ and $T_x M$. The horisontal subspaces $H(x, \Lambda)$ can be seen as images of linear maps $h(x, \Lambda) \colon T_x M \to T_{(x,\Lambda)} D^{r,l}(M, M')$ such that $d\pi \circ h(x, \Lambda) = \mathsf{id}$ for all $(x, \Lambda) \in D^{r,l}(M, M')$. It is easy to see that the correspondence between $H(x, \Lambda)$ and $h(x, \Lambda)$ is one-to-one and that the maps $h(x, \Lambda)$ as above can be seen as sections of the associated affine bundle over $D^{r,l}(M, M')$ (this is a general property: all connections on a fibration can be seen as sections in an affine bundle over the fibration space).

By Theorem 2.1, the $r$-jet section $s(x) := j_x^k f$ is horisontal for any $l$-nondegenerate CR-map $f \colon U \to M'$, where $U \subset M$ is any open subset. Moreover, $s(x) = (\Lambda_i^\alpha(x))$ always satisfies the contact system $\partial_{x_j} \Lambda_i^\alpha(x) = \Lambda_i^{\alpha + e_j}(x)$, where $e_j \in \mathbb{Z}_+^{\dim M}$ is the unit vector whose $j$th component is 1 and other components are 0. We now keep the highest terms $h_i^\alpha$ for $|\alpha| = r$ unchanged and replace the lower terms by the contact system. Then the section $s$ still satisfies the new system, also called $h$ and any other horisontal section of $D^{r,l}(M, M')$ with respect to $h$ over $U$ is an $r$-jet section of some map $f \colon U \to M'$.

In order to obtain the conclusion of Corollary 2.4 we change the connection given by $h$ as follows. Let $S \subset D^{r,l}(M, M')$ be the submanifold of all $r$-jets whose first order parts satisfy the Cauchy-Riemann equations on $M$, i.e. send $T^c M$ into $T^c M'$ and are $\mathbb{C}$-linear there. Clearly, the $r$-jet lift $\{(x, j_x^r f) : x \in U\}$ of any CR-map $f \colon U \to M'$ as above is contained in $S$. By using a partition of unity we can define a real smooth function $g \colon D^{r,l}(M, M') \to [0, 1] \subset \mathbb{R}$ such that $S = \{g = 0\}$.

Let $s \colon U \to D^{r,l}(M, M')$ be any horisontal section, i.e. $H(x, s(x))$ coincides with the graph of $ds(x)$ for every $x \in U$. In local coordinates, this condition can be written as $\partial_{x_j} s(x) = h_j(x, s(x))$. The symmetry property $\partial_{x_i} \partial_{x_j} s = \partial_{x_j} \partial_{x_i} s$ implies now the integrability condition $\partial_{x_i} h_j(x, s(x)) =$

$\partial_{x_j} h_i(x, s(x))$ or

$$\partial_{x_i} h_j + (\partial_\Lambda h_j) h_i = \partial_{x_j} h_i + (\partial_\Lambda h_i) h_j, \tag{5.2}$$

where we used the matrix notation and substituted $h_i$ and $h_j$ for $\partial_{x_i} s$ and $\partial_{x_j} s$ respectively. We now write $h_0$ for the system given by the map $\Phi$ as above and put $h := (1 - g)h_0 + gh_1$, where $h_1$ is another connection to be chosen and $g$ is the function defined above. Due to the affine structure on the space of connections, $h$ is well-defined. Then, since dim $M \geq 3$, the Thom transversality theorem (see e.g. [Hi76]) yields the existence of $h_1$ such that (5.2) does not hold outside $S$. On the other hand, we have $h = h_0$ on $S$ by the construction. Hence, if $s(x)$ is horisontal with respect to the connection given by $h$, its graph must be contained in $S$ and moreover $s$ is the $r$-jet section of some map $f \colon U \to M'$ by the construction of $h$. Finally, by the construction of $S$, $f$ is CR as required. $\qquad\square$

Mathematisches Institut, Eberhard-Karls-Universität Tübingen, 72076 Tübingen, Germany
*E-mail address*: sykim@math.snu.ac.kr, dmitri.zaitsev@uni-tuebingen.de